\documentclass[preprint,12pt]{article}

\usepackage{graphicx}%
\usepackage{multirow}%
\usepackage{amsmath,amssymb,amsfonts}%
\usepackage{amsthm}%
\usepackage{mathrsfs}%
\usepackage[title]{appendix}%
\usepackage[table]{xcolor}%
\usepackage{textcomp}%
\usepackage{manyfoot}%
\usepackage{booktabs}%
\usepackage{algorithm}%
\usepackage{algorithmicx}%
\usepackage{algpseudocode}%
\usepackage{listings}%
\usepackage{algpseudocode}
\usepackage{svg}
\usepackage{subcaption}    
\usepackage{a4wide}
\usepackage{booktabs}
\usepackage{tikz-cd}
\usepackage{url}

\theoremstyle{definition}%
\newtheorem{theorem}{Theorem}
\newtheorem{proposition}[theorem]{Proposition}%
\newtheorem{corollary}[theorem]{Corollary}
\newtheorem{lemma}[theorem]{Lemma}

\theoremstyle{definition}%
\newtheorem{example}{Example}%

\theoremstyle{definition}%
\newtheorem{definition}{Definition}%

\newcommand{\R}{\mathbb{R}}
\newcommand{\Ess}{\operatorname{Ess}}
\newcommand{\Ord}{\operatorname{Ord}}
\newcommand{\Rel}{\operatorname{Rel}}

\newcommand{\ord}{\operatorname{ord}}
\newcommand{\rel}{\operatorname{rel}}
\newcommand{\XPHT}{\operatorname{XPHT}}
\newcommand{\PHT}{\operatorname{PHT}}

\title{An Implementation of the 3-Dimensional Extended Persistent Homology Transform}

\author{
  Kiyoshi Takeuchi Romo, Vanessa Robins \\
  Research School of Physics \\
  The Australian National University \\
  Canberra\\
  \texttt{\{kiyoshi.takeuchiromo, vanessa.robins\}@anu.edu.au} \\
}

\begin{document}

\maketitle

\begin{abstract}The Extended Persistent Homology Transform (\(\mathrm{XPHT}\)) describes a shape through the extended persistence diagrams of its directional height functions. We present an algorithm for computing a finite directional sample of the three-dimensional \(\mathrm{XPHT}\) from a triangulated boundary. Our main construction determines essential first-homology pairings by combining an ordered tree–cotree decomposition with evaluation-dual cocycles. The decomposition supplies a homology basis compatible with the height filtration, while the cocycles are used to determine the essential birth–death pairs. Degree-zero persistence and duality determine the remaining boundary diagrams. Under the boundary-reduction hypotheses and with critical-point sign labels, these diagrams yield the extended persistence of the solid. The resulting procedure uses graph constructions and linear algebra over $\mathbb F_2$ on boundary meshes extracted from three-dimensional binary images.
\end{abstract}

\section{Introduction}

Persistent homology is a family of tools used in topological data analysis to track how the homology of a space $X$ changes as a parameter varies. These changes are captured in what is called a filtration of the space $X$. A filtration, in its most basic sense, is a sequence of nested subspaces $X_0\subseteq \cdots\subseteq X$. The information we obtain varies depending on which filtration we choose. Common examples include the Vietoris–Rips filtration and the height filtration. The Vietoris–Rips filtration is frequently used to infer the underlying topological structure of point clouds \cite{Rabadan_Blumberg_2019}, whereas a height filtration examines the sublevel sets of a height function on a given space \cite{edelsbrunner2010computational}. The latter depends on the direction chosen to represent ``upwards.''

In standard Persistent Homology, we have sequences of the form \[
H_k(X_{a_0})\to \cdots \to H_k(X_{a_m})
\]
where $H_k$ is the $k$-th homology space. The critical points of the height function mark the levels at which homology classes can appear or disappear in the height filtration. We then pair some of the critical points of the height function in a meaningful algebraic and geometric sense. Using this sequence usually leaves some critical points of the height function unpaired.
On the other hand, in Extended Persistent Homology we have ``extended'' sequences of the form
\[
H_k(X_{a_0})\to \cdots \to H_k(X_{a_m})
\to
H_k(X,X^{a_m})\to \cdots \to H_k(X,X^{a_0}).
\]
Given a Morse height filtration on a compact finitely triangulated $3$-manifold with boundary, one may pair all critical points \cite{cohensteiner2009extending}.

The $\XPHT$ is a function that packages extended persistent homology information obtained from height filtrations in all possible directions. Our focus is on implementing an algorithm to compute a discrete version of the $3$D $\XPHT$ of a digital image. As shown in \cite{Turner2014PHT}, the $\PHT$ (and consequently the $\XPHT$) provides a faithful summary of a shape within the class considered, since the original shape can be recovered from its transform.

The $2$D $\XPHT$ was implemented in \cite{Turner2024XPHT} as a computer algorithm using a duality theorem that substantially simplifies the calculation: only zeroth homology needs to be computed explicitly. In the $3$D case, the duality theorem still provides considerable information, but an additional step is required to compute the first homology classes and determine the pairing of essential classes. We address this problem in Section~\ref{chapter-pairing}. 

The structure of this paper goes as follows; We start by giving an overview of works that are relevant in Section \ref{chapter-previouswork}. In Section~\ref{chapter-background}, we introduce the algebraic foundations of extended persistence, the evaluation pairing between homology and cohomology, and the duality relations used to reconstruct the boundary diagrams. We provide all the necessary definitions and algorithms to compute the discretized version of the $\XPHT$ starting from a $3$D binary image. 

In Section \ref{chapter-xpht0}, we tackle the part of the algorithm corresponding to $\XPHT_0$. We spend some amount of time here, since it will also be the foundation for the  algorithm for the essential $1$ classes. As we prove in Corollary \ref{simplification}, $\XPHT_0$ alone provides a wealth of information for our $3$D shape. In Section \ref{chapter-h1classes}, we implement the tree-cotree algorithm to obtain the essential $1$ classes. In Section \ref{chapter-pairing}, we have our main theoretical contribution: the pairing of the essential $1$ classes. 

Essential $1$ classes on surfaces can be paired using
the Reeb graph of the height function, as reviewed in
\cite[Section~2]{cohensteiner2009extending}. For the sampled
$\XPHT$, this approach requires obtaining the corresponding
Reeb graph for each directional height function. Our method computes the essential first-homology pairing
directly from the triangulated boundary. The resulting algorithm uses graph constructions and linear algebra on the boundary mesh,
avoiding the construction of Reeb graphs.

We end the manuscript with conclusions and future work, and note that a working implementation of the discretized algorithm described in this manuscript is available in \cite{kiyoshiandres_3dxpht}. 

\section{Previous Work}\label{chapter-previouswork}

Persistent homology has been applied across a wide range of scientific disciplines, including medicine \cite{BritoPacheco2026}, astronomy \cite{XU201934}, and geology \cite{https://doi.org/10.1002/2015WR017937}. These applications illustrate the value of topological summaries for characterizing complex data. In shape analysis, the $\XPHT$ combines such summaries across directional height filtrations, retaining information about both topology and geometry. Its application to detecting planar symmetries \cite{inproceedings} demonstrates the utility of this directional information and motivates the development of computational methods for three-dimensional shapes.

General-purpose software already provides tools for computing extended persistent homology. For example, GUDHI \cite{gudhi} supports extended persistence for lower-star filtrations induced by values assigned to the vertices of a simplicial complex. This framework includes height filtrations, obtained by assigning each vertex its height in a chosen direction. Our work addresses the more specific task of computing the $3$D $\XPHT$ by exploiting the structure of the underlying shape and resolving the pairing of essential first homology classes.

To address this pairing problem, we build on the tree–cotree framework described in \cite{Eppstein2003} and developed further in \cite{Erickson2012Cycles}. This framework provides a combinatorial approach to constructing cycle representatives that form a basis for first homology. Obtaining such a basis, however, does not by itself determine the birth–death pairings required for extended persistence. Our contribution builds on these constructions to determine the pairing of essential $1$ homology classes, as described in Section~\ref{chapter-pairing}.

 \cite{EdelsbrunnerOlsbock2020} generalize tree--cotree decompositions to ordered polyhedral complexes over $\mathbb{F}_2$. On a closed triangulated surface, their degree-one construction specializes, with matching ordering conventions, to the ordered tree--cotree basis used here. They construct homology and cohomology bases dual under evaluation using exhaustive matrix reduction, whereas we obtain evaluation-dual cocycle representatives for our chosen homology basis by solving linear systems.

A related computational perspective comes from persistent cohomology. The relationships between persistent homology and cohomology were developed in \cite{Silva2011}, establishing how equivalent persistence information can be obtained through these different algebraic formulations. This provides a foundation for using cohomological methods in persistence computations and connects our approach to a broader family of methods that exploit duality.

A computational implementation also requires restricting the calculation to a finite set of directions on the sphere. Sufficient conditions on the number of directions needed to determine a shape are discussed in \cite{curry2021directionsdetermineshapesufficiency}. For $3$D simplicial complexes and a variation of persistent homology, it is discussed in \cite{fasy_faithful_2026}. We further discuss this in \ref{sampling-directions}.

For the purpose of this work, we assume that there is an algorithm to obtain a bounding surface of the $3$D manifold that is fairly accurate with respect to the positions of  critical points. A popular visualization method, called marching cubes \cite{wenger2013isosurfaces}, is one such algorithm through which an approximation to the boundary may be obtained of the $3$-manifold. We note that this method does not guarantee that the topology of this approximation will be the same as that of the boundary of the manifold. In particular one may end up with different connectivity among other differences.

\section{Background}\label{chapter-background}

The Persistent Homology Transform ($\PHT$) was introduced in \cite{Turner2014PHT}, where it was shown to be injective on the class of shapes considered. Consequently, a shape in this class is uniquely determined by its collection of persistence diagrams indexed by direction.

\subsection{Foundations}

The $\XPHT$ can be defined in a general setting using the framework of definable sets \cite{Coste_1999}. However, for the discrete implementation developed here, it suffices to work with finite triangulated spaces. We therefore begin by fixing the notation for the underlying manifold, its boundary, and a triangulation of that boundary.

Let $X\subset\mathbb{R}^3$ be a compact embedded
$3$-manifold with boundary $A=\partial X$. Let $K$ be a finite simplicial complex triangulating $A$, so that its geometric realization $|K|$ is homeomorphic to $A$. This triangulation provides a combinatorial representation of the boundary on which homology can be computed.

With coefficients in a field, the homology groups of $K$ are vector spaces. Throughout this paper, we use coefficients in $\mathbb{F}_2$, which simplifies computation by eliminating the need to track orientation signs in boundary calculations. Unless the coefficient field requires emphasis, we omit it from the notation and write $H_i(K)$ for $H_i(K;\mathbb{F}_2)$.

\subsection{Homology and Cohomology}

We recall the definitions and duality results needed for our construction. For a more comprehensive treatment of homology and cohomology, we refer to \cite{MR1867354}. All homology and cohomology groups are taken with coefficients in $\mathbb{F}_2$, unless stated otherwise.

An \textbf{abstract simplicial complex} on a vertex set $V$ is a collection $K$ of finite subsets of $V$ that is closed under taking subsets: if $\sigma\in K$ and $\tau\subseteq\sigma$, then $\tau\in K$. A simplex $\sigma$ with $|\sigma|=k+1$ is called a \textbf{$k$-simplex}. We work with finite simplicial complexes triangulating manifolds, together with the subcomplexes arising in their filtrations.

The vector space $C_k(K;\mathbb{F}_2)$ generated by the $k$-simplices is the space of \textbf{$k$-chains}. The \textbf{boundary operator} is the linear map
\begin{align*}
\partial_k:C_k(K)&\longrightarrow C_{k-1}(K),\\
[v_0,\ldots,v_k]&\longmapsto
\sum_{r=0}^{k}[v_0,\ldots,\widehat{v_r},\ldots,v_k],
\end{align*}
where the hat indicates an omitted vertex. We use ordinary, unreduced homology, so $\partial_0=0$. Since $\partial_k\partial_{k+1}=0$, every boundary is a cycle. The \textbf{$k$-th homology group} is therefore the vector space
\[
H_k(K)=\frac{\ker\partial_k}{\operatorname{im}\partial_{k+1}}.
\]

The dual vector space
\[
C^k(K)=\operatorname{Hom}(C_k(K),\mathbb{F}_2)
\]
is the space of \textbf{$k$-cochains}. Its \textbf{coboundary operator} is
\begin{align*}
\partial^k:C^k(K)&\longrightarrow C^{k+1}(K),\\
\alpha&\longmapsto\alpha\circ\partial_{k+1}.
\end{align*}
The identity $\partial^k\partial^{k-1}=0$ gives the \textbf{$k$-th cohomology group}
\[
H^k(K)=\frac{\ker\partial^k}{\operatorname{im}\partial^{k-1}}.
\]
We identify these groups with those of the triangulated space when convenient.

Over a field, cohomology is naturally identified with the dual of homology. This identification provides the evaluation-dual cocycles used in Section~\ref{chapter-pairing}.

\begin{theorem}\label{uci}
For a finite simplicial complex $K$, the evaluation map
\begin{align*}
UC:H^1(K)&\longrightarrow\operatorname{Hom}(H_1(K),\mathbb{F}_2),\\
[\alpha]&\longmapsto\bigl([z]\longmapsto\alpha(z)\bigr)
\end{align*}
is an isomorphism.
\end{theorem}

\begin{proof}
\cite[p.~198]{MR1867354}.
\end{proof}

\begin{corollary}\label{evaluation-dual-basis}
Let $z_1,\ldots,z_r$ be cycles whose classes form a basis of $H_1(K)$. There exist cocycles $\widehat z_1,\ldots,\widehat z_r$ satisfying
\[
\widehat z_i(z_j)=\delta_{ij}.
\]
Their cohomology classes form the basis of $H^1(K)$ dual to $[z_1],\ldots,[z_r]$ under the evaluation pairing.
\end{corollary}

\begin{proof}
Apply Theorem~\ref{uci} to the algebraic dual basis of $\operatorname{Hom}(H_1(K),\mathbb{F}_2)$, and choose a cocycle representative for each resulting cohomology class.
\end{proof}

\subsection{Persistent Homology}

We now describe the filtrations used to construct the transforms. To distinguish a general input shape from the $3$-manifold $X$ and its boundary $A$, let $Y\subset\mathbb{R}^n$ denote a compact shape represented by a finite geometric simplicial complex. Here $n$ is the ambient dimension, not necessarily the intrinsic dimension of $Y$.

For a direction $u\in S^{n-1}$, define the height function
\[
h_u:Y\longrightarrow\mathbb{R},\qquad h_u(x)=x\cdot u.
\]
Its \textbf{sublevel sets} and \textbf{superlevel sets} are
\[
Y_t=h_u^{-1}(( -\infty,t]),\qquad
Y^t=h_u^{-1}([t,\infty)).
\]
For $s\le t$, inclusion induces a linear map $H_k(Y_s)\to H_k(Y_t)$. These vector spaces and maps form the \textbf{$k$-th persistence module}
\[
\mathsf{PH}_k(Y,h_u)=\bigl(H_k(Y_t)\bigr)_{t\in\mathbb{R}}.
\]

In this finite piecewise-linear setting, the module changes only at finitely many parameter values and admits an interval decomposition. The general interval decomposition theorem is given in \cite{crawleyboevey2014decompositionpointwisefinitedimensionalpersistence}; finiteness in the definable setting is also discussed in \cite{curry2021directionsdetermineshapesufficiency}. The resulting \textbf{persistence diagram} $D_k(Y,h_u)$ is a multiset of birth--death pairs $(b,d)$, counted with multiplicity. For the sublevel-set convention used here, finite bars have the form $[b,d)$, while bars that persist indefinitely have the form $[b,\infty)$. The latter represent \textbf{essential classes}.

The \textbf{$k$-th Persistent Homology Transform} \cite{Turner2014PHT} is the map
\begin{align*}
\mathrm{PHT}_k(Y):S^{n-1}&\longrightarrow\mathcal{D},\\
u&\longmapsto D_k(Y,h_u),
\end{align*}
where $\mathcal{D}$ is the space of persistence diagrams. Packaging degrees $0,\ldots,n-1$ gives
\begin{align*}
\mathrm{PHT}(Y):S^{n-1}&\longrightarrow\mathcal{D}^{n},\\
u&\longmapsto\bigl(D_0(Y,h_u),\ldots,D_{n-1}(Y,h_u)\bigr).
\end{align*}

A \textbf{critical value} is a parameter at which the topology of the sublevel filtration changes. Since $h_u$ is linear on each simplex, such changes can occur only at vertex heights, although not every vertex height need be critical. Equal heights may be processed using a consistent ordering of simplices that places faces before cofaces. Such an ordering refines the computation; it does not, by itself, change the assigned height values.

Let
\[
\mathbb R_{\operatorname{ext}}
=
(\mathbb R\times\{\ord\})\cup(\mathbb R\times\{\rel\})
\]
with the following total ordering:
\[
(s,\ord)\le (t,\ord)\quad\text{if }s\le t,
\qquad
(s,\ord)\le (t,\rel)\quad\text{for all }s,t,
\]
and
\[
(s,\rel)\le (t,\rel)\quad\text{if }s\ge t.
\]
Thus the ordinary part of the filtration is read in increasing height, while the relative part is read in decreasing height.

The \textbf{$k$-th extended persistent homology} of $(X,h_u)$ is defined by the sequence as in \cite{cohensteiner2009extending}:
\[
H_k(X_{a_0})\to \cdots \to H_k(X_{a_m})
\to
H_k(X,X^{a_m})\to \cdots \to H_k(X,X^{a_0}),
\]
where $a_0<\cdots<a_m$ are filtration values containing all critical values. The first half is called the \textbf{ordinary part}, and the second half is called the \textbf{relative part}.

The intervals in the decomposition of this persistence module come in three types. An \textbf{ordinary interval} is born and dies in the ordinary part:
\[
[(b,\ord),(d,\ord)).
\]
An \textbf{essential interval} is born in the ordinary part and dies in the relative part:
\[
[(b,\ord),(d,\rel)).
\]
A \textbf{relative interval} is born and dies in the relative part:
\[
[(b,\rel),(d,\rel)).
\]
We denote the corresponding collections of intervals by
\[
\Ord_k(X,h_u),\qquad \Ess_k(X,h_u),\qquad \Rel_k(X,h_u),
\]
and omit $h_u$ when no confusion arises.
Analogously we have the \textbf{$k$-th Extended Persistent Homology Transform} $(\XPHT)$, which is the map
\begin{align*}
    \XPHT_k(X):S^{n-1}&\longrightarrow \mathcal D_{\operatorname{ext}}\\
    u&\longmapsto
    \bigl(\Ord_k(X,h_u),\Ess_k(X,h_u),\Rel_k(X,h_u)\bigr),
\end{align*}
where $\mathcal D_{\operatorname{ext}}$ denotes the space of extended persistence diagrams.

Often one packages all dimensions together:
\begin{align*}
    \XPHT(X):S^{n-1}&\longrightarrow \mathcal D_{\operatorname{ext}}^n\\
    u&\longmapsto
    \bigl\{\XPHT_0(X)(u),\ldots,\XPHT_{n-1}(X)(u)\bigr\}.
\end{align*}

The following results relate homological degrees and opposite height directions. We state the manifold symmetry results in the closed-manifold setting needed for the boundary surface $A$. The manifold dimension is denoted by $m$, reserving $b$ and $d$ for interval endpoints. All correspondences preserve multiplicity.

\begin{theorem}[Symmetry; Proposition~2.8 of \cite{Turner2024XPHT}]\label{duality-1}
Let $M$ be a closed triangulated $m$-manifold and let $f:M\to\mathbb{R}$ be piecewise linear. There are bijections of interval multisets
\begin{align*}
\Ord_k(M,f)&\longleftrightarrow\Ord_{m-k-1}(M,-f),\\
\Ess_k(M,f)&\longleftrightarrow\Ess_{m-k}(M,-f),\\
\Rel_k(M,f)&\longleftrightarrow\Rel_{m-k+1}(M,-f).
\end{align*}

In $\Ord$ and $\Rel$, a birth--death pair $(b,d)$ is sent to $(-d,-b)$. For $\Ess$, $[(b,\ord),(d,\rel) )$ is sent to $[(-d,\ord),(-b,\rel))$.
\end{theorem}

\begin{theorem}[Ordinary--relative correspondence; Corollary~2.7 of \cite{Turner2024XPHT}]\label{duality-2}
Let $M$ be a finite simplicial complex and let $f:|M|\to\mathbb{R}$ be linear on each simplex. There is a bijection
\[
\Ord_k(M,f)\longleftrightarrow\Rel_{k+1}(M,-f)
\]
that sends an ordinary birth--death pair $(b,d)$ to the relative pair $(-b,-d)$.
\end{theorem}

\begin{theorem}[Complementary-degree correspondence; proof of Proposition~2.8 of \cite{Turner2024XPHT}]\label{duality-3}
Under the closed-manifold assumptions of Theorem~\ref{duality-1}, there is a bijection
\[
\Ord_{m-k}(M,f)\longleftrightarrow\Rel_k(M,f)
\]
that sends an ordinary birth--death pair $(b,d)$ to the relative pair $(d,b)$.
\end{theorem}

For height functions, changing $u$ to $-u$ changes $h_u$ to $-h_u$, so these results allow information to be transferred between opposite directions. When $X \subset \R^2$, its boundary is one-dimensional, and boundary duality reduces the persistence calculation to degree zero. In three dimensions, the boundary is a surface, and essential first homology remains to be computed and paired.


The following results are simple consequences of the general results established in \cite{Turner2024XPHT}, and apply only to our $3$D setting.

\begin{lemma}\label{ord-2}
For a closed triangulated surface $A$, the multiset $\Ord_2(A,h_u)$ is empty for every direction $u$.
\end{lemma}

\begin{proof}
The sublevel filtration can be computed using the lower-star subcomplexes of a triangulation. Since the complex is two-dimensional, it contains no $3$-chains. Thus a nonzero $2$-cycle cannot become a boundary at a later stage. All maps on second homology induced by these subcomplex inclusions are injective, so no second homology class dies in the ordinary part of the filtration.
\end{proof}

\begin{corollary}\label{simplification}
Let $A=\partial X$ be the compact, oriented boundary surface of a $3$-manifold $X \subset \R^3$. Computing $\Ord_0(A,h_u)$, $\Ess_0(A,h_u)$, and $\Ess_1(A,h_u)$ for all directions $u$ is sufficient to recover the full transform $\XPHT(A)$. For a finite directional sample, the same conclusion holds provided that the sample also contains the required antipodal directions.
\end{corollary}

\begin{proof}
Only degrees $0$, $1$, and $2$ occur. Apply Theorem~\ref{duality-1} with $m=2$ to recover $\Ord_1(A,h_u)$ from $\Ord_0(A,h_{-u})$. Lemma~\ref{ord-2} gives $\Ord_2(A,h_u)=\varnothing$. Theorems~\ref{duality-2} and~\ref{duality-3} determine the relative intervals from these ordinary intervals. Finally, essential symmetry recovers $\Ess_2(A,h_u)$ from $\Ess_0(A,h_{-u})$. The remaining essential intervals are precisely those in degree one.
\end{proof}

The resulting reconstruction rules are summarized in Table~\ref{tab:xpht}. In each nonempty row, $(b,d)$ denotes an interval in the source multiset.

\begin{table}[htbp]
\centering
\begin{tabular}{lll}
\hline
Target at $u$ & Source & Endpoint rule\\
\hline
$\Ord_0(A,h_u)$ & Computed directly & $(b,d)$\\
$\Ord_1(A,h_u)$ & $\Ord_0(A,h_{-u})$ & $(-d,-b)$\\
$\Ord_2(A,h_u)$ & Empty & ---\\
$\Ess_0(A,h_u)$ & Computed directly & $(b,d)$\\
$\Ess_1(A,h_u)$ & Computed directly & $(b,d)$\\
$\Ess_2(A,h_u)$ & $\Ess_0(A,h_{-u})$ & $(-b,-d)$\\
$\Rel_0(A,h_u)$ & Empty & ---\\
$\Rel_1(A,h_u)$ & $\Ord_0(A,h_{-u})$ & $(-b,-d)$\\
$\Rel_2(A,h_u)$ & $\Ord_0(A,h_u)$ & $(d,b)$\\
\hline
\end{tabular}
\caption{Recovery of the extended persistence of the boundary surface from degree-zero data and essential degree-one data.}
\label{tab:xpht}
\end{table}

Assume that $h_v$ satisfies the hypotheses of the
boundary-reduction theorem
\cite[Theorem~4.18]{Turner2024XPHT}.
Let $I$ be an ordinary or essential degree-$k$ interval
of the boundary $A=\partial X$, born at a boundary critical point $p$ for $h_v$ on the solid. Write
$\eta(p)\in\{-1,+1\}$ for the sign defined by the
corresponding boundary Morse chart.
The boundary-reduction theorem identifies $I$ as an
interval belonging to the solid $X$ precisely when
$\eta(p)=+1$.

In the smooth setting, this sign has a direct
interpretation using the outward unit normal $n(p)$.
At a boundary critical point, the derivative of the
height function in the inward normal direction is
\[
D_{-n(p)}h_v=-n(p)\cdot v.
\]
Consequently, $p$ is positive critical precisely when
$n(p)\cdot v<0$.

For the triangulated boundary, we estimate this sign
using an averaged outward normal. Let $\mathcal F(p)$
denote the triangular faces incident to $p$, and let
$n_\tau$ be the outward unit normal of face $\tau$.
Using equal weights, define
\[
\bar n(p)
=
\frac{1}{|\mathcal F(p)|}
\sum_{\tau\in\mathcal F(p)}n_\tau.
\]
All face normals are oriented outward from the solid
$X$; on a cavity boundary, they therefore point into
the cavity. Whenever $\bar n(p)\cdot v\ne0$, the
estimated critical-point sign is
\[
\widehat\eta(p)
=
-\operatorname{sgn}\bigl(\bar n(p)\cdot v\bigr).
\]

The normal-based selection is made under the
assumption that these estimates agree with the
boundary Morse signs at every ordinary birth event
used in the selection:
\[
\widehat\eta(p)=\eta(p).
\]
In the case of essential $1$ cycles, we instead analogously define the sign of an edge. We may assume that $\bar n(e)\neq 0$ because, in the case of $\widehat \eta(e)\neq 0$ for the leftover edge $e$ (c.f. Section \ref{chapter-h1classes}), we may replace the edge with one of equal height sharing a vertex, belonging to the fundamental tree, belonging to the same cycle, and having nonzero sign.

Under this assumption, retaining exactly those
ordinary and essential boundary intervals whose
birth vertices satisfy
\[
\bar n(p)\cdot v<0
\]
implements the selection prescribed by the
boundary-reduction theorem.

The relative intervals of $X$ are recovered by
applying Theorem \ref{duality-2} to its selected ordinary
intervals in the antipodal direction. For the purpose of the implementations, assume that the vertex sign is transported throughout the algorithms.

\subsection{Computational Background}

Our implementation combines connectivity tracking, surface triangulation, directional sampling, and linear algebra over $\mathbb{F}_2$. We first introduce the union--find data structure, which supports the construction of spanning trees and cotrees. We then describe the representation of digital images and the choice of sampling directions, before explaining how the resulting mesh is used to construct homology bases and boundary matrices.

\subsubsection{Union--Find}

The \textbf{union--find} data structure, also called a disjoint-set data structure, maintains a partition of a collection of elements into disjoint subsets \cite{sedgewick1988algorithms}. It supports two operations: $\textsc{Find}(x)$ returns a representative of the subset containing $x$, and $\textsc{Union}(x,y)$ merges the subsets containing $x$ and $y$. When processing a graph, these subsets represent the connected components formed by the edges processed so far.

The subsets are represented by rooted trees. Each element $x$ stores a pointer $\mathrm{parent}[x]$ to its parent, and a root satisfies $\mathrm{parent}[x]=x$. Following parent pointers therefore leads to the representative of the component. To accelerate subsequent queries, \textsc{Find} uses \textbf{path compression}, updating the pointers along a search path to point directly to the root.

The \textsc{Union} operation uses \textbf{union by rank}. The quantity $\mathrm{rank}[x]$ is a balancing quantity associated with a root; it is not a filtration rank or an ordering of the vertices. When two roots have different ranks, the root of smaller rank is attached to the root of larger rank. When their ranks agree, either root may be retained, and its rank is increased by one.

In persistence computations, the representative chosen for balancing need not be the vertex at which a component was born. Birth information must therefore be maintained separately when applying the elder rule.

\subsubsection{Digital Images and Boundary Meshes}

Following \cite{robins2011discrete}, a \textbf{three-dimensional digital image} is a function $f:D\to\mathbb{R}$ defined on a finite rectangular grid

$$
D=
\left\{
(i,j,k)\in\mathbb{Z}^3:
0\le i<N_x,\quad
0\le j<N_y,\quad
0\le k<N_z
\right\}.
$$

Each grid index identifies a voxel. The image is \textbf{binary} when $f(D)\subseteq\{0,1\}$, with $1$ representing the foreground and $0$ the background.

To apply the surface-based constructions developed above, we extract a triangular mesh representing the interface between foreground and background. Marching-cubes methods provide a standard approach to this task \cite{wenger2013isosurfaces}. For binary data, an isovalue of $1/2$ separates the two voxel values. If voxel spacings are unequal, the mesh coordinates must incorporate those spacings so that subsequent height functions are evaluated in physical coordinates.

Throughout the following constructions, we assume that the extracted mesh $K$ is an embedded, closed, orientable triangulated surface representing $A=\partial X$. In particular, every edge is incident to exactly two triangular faces, and the link of every vertex is a circle. These are requirements on the mesh used by the algorithm; they should not be assumed for every marching-cubes implementation or every input configuration.

To avoid truncating the extracted surface at the boundary of the image domain, the foreground should be surrounded by background voxels, adding padding when necessary.

\subsubsection{Sampling Directions}\label{sampling-directions}

A computational implementation evaluates the $\XPHT$ on a finite set of directions. Motivated by spherical sampling constructions such as \cite{doi:10.1073/pnas.2426574122}, we obtain a symmetric directional sample by subdividing a regular octahedron and projecting its vertices radially onto the unit sphere.

\begin{figure}[htbp]
\centering
\includegraphics[width=0.5\linewidth]{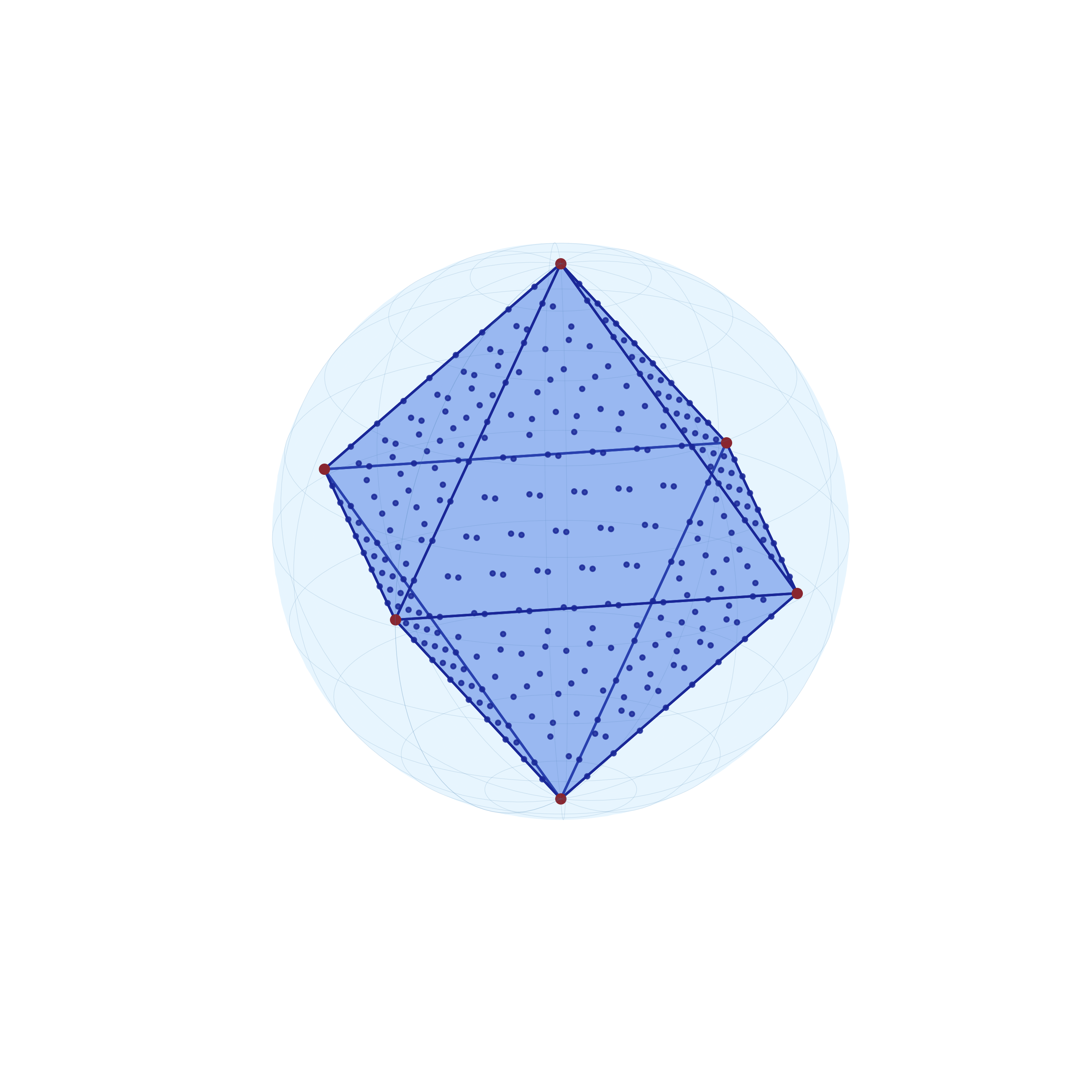}
\caption{A set of $326$ directions obtained by radially projecting the vertices of a subdivided regular octahedron onto the unit sphere. Each original edge is divided into nine equal segments.}
\label{fig:octahedron}
\end{figure}

Let $q\ge1$ be an integer. Subdivide each edge of the octahedron into $q$ equal segments and each triangular face into the corresponding triangular grid \ref{fig:octahedron}. The resulting vertices lie in

$$
V_q=
\left\{
\frac{1}{q}(i,j,k):
(i,j,k)\in\mathbb{Z}^3,\quad
|i|+|j|+|k|=q
\right\}.
$$

Their radial projections define the direction set

$$
U_q=
\left\{
\frac{(i,j,k)}{\sqrt{i^2+j^2+k^2}}:
(i,j,k)\in\mathbb{Z}^3,\quad
|i|+|j|+|k|=q
\right\}
\subset S^2.
$$

This set is antipodally symmetric: if $u\in U_q$, then $-u\in U_q$. Thus the opposite directions needed by the duality relations are included automatically.

For $q=9$, this gives $326$ directions. Radial projection gives a symmetric covering of the sphere, although the resulting spherical cells do not all have equal area. See \cite{monteverde2024platonic} for an implementation in MATLAB.

\subsubsection{Linear Algebra}

The chain and cochain operations introduced above can be represented by matrices over $\mathbb{F}_2$. Fix an ordering of the simplices in each dimension, and use the corresponding simplex bases for chains and dual bases for cochains.

\begin{definition}
Let $\tau_1,\ldots,\tau_m$ be the $(k-1)$-simplices of $K$, and let $\sigma_1,\ldots,\sigma_n$ be its $k$-simplices. The \textbf{$k$-th boundary matrix} is the matrix $B_k\in\mathbb{F}_2^{m\times n}$ with entries

$$
(B_k)_{ij}
=
\begin{cases}
1, & \tau_i\text{ is a codimension-one face of }\sigma_j,\\
0, & \text{otherwise}.
\end{cases}
$$

The \textbf{$(k-1)$-th coboundary matrix}, denoted $B^{k-1}$, represents $\partial^{k-1}:C^{k-1}(K)\to C^k(K)$ in the dual simplex bases.
\end{definition}

In particular, $B_1$ has vertices as rows and edges as columns, while $B_2$ has edges as rows and triangular faces as columns.

\begin{proposition}\label{boundary_matrix}
Let $\operatorname{coord}(c)$ denote the coordinate column vector of a chain or cochain in its chosen basis. Then

$$
\operatorname{coord}(\partial_k c)
=
B_k\operatorname{coord}(c),
\qquad
B^{k-1}=B_k^\top.
$$

Consequently,

$$
\operatorname{coord}(\partial^{k-1}\alpha)
=
B_k^\top\operatorname{coord}(\alpha).
$$

\end{proposition}

\begin{proof}
See \cite{edelsbrunner2010computational} for a proof.
\end{proof}

\section{$\XPHT_0$}\label{chapter-xpht0}

Fix a direction $v\in S^2$ and consider the height function
$h_v(x)=x\cdot v$ on the triangulated boundary $K$.
We compute degree-zero persistence by tracking the connected
components of its sublevel sets. Only vertices and edges are
needed for this computation, since adding higher-dimensional
simplices does not change connectivity.

\subsection{Tie-Breaking Convention}\label{total-order}

Assign each vertex $p$ the filtration value $h_v(p)$ and
each edge $e=[p,q]$ the value

$$
h_v(e)=\max\{h_v(p),h_v(q)\}.
$$

Thus an edge enters the filtration when both of its endpoints
have entered. To make the computation deterministic, we fix a
\textbf{tie-breaking convention}. Vertices are ordered first
by height and then lexicographically by their coordinates:
$$
\kappa(p)=\bigl(h_v(p),p_x,p_y,p_z\bigr).
$$

We assume that no two vertices share the same coordinates, as this may lead to degenerate situations. At each filtration level $t$, all vertices at height $t$ are inserted in the vertex tie-breaking order before any edges entering at that level. For an edge $e=[p,q]$, let $i(e)<j(e)$ be the fixed indices of its endpoints, and define
\[
\varepsilon(e)=
\begin{cases}
0, & h_v(p)=h_v(q),\\
1, & h_v(p)\ne h_v(q).
\end{cases}
\]
We order edges lexicographically by
\[
\lambda(e)
=
\bigl(h_v(e),\varepsilon(e),i(e),j(e)\bigr).
\]
Thus, among edges entering at the same height, edges
whose endpoints have equal heights are processed first.
Within each group, edges are ordered lexicographically
by their endpoint indices. 

\subsection{Tracking Births and Mergers}

For the closed boundary surface considered here,
$\Rel_0(K,h_v)=\emptyset$ by
Corollary~\ref{simplification}. It therefore remains to compute
$\Ord_0(K,h_v)$ and $\Ess_0(K,h_v)$.

We use a union--find structure to maintain the connected
components of the current sublevel set. Each component stores
a birth vertex and its birth height. These are maintained
separately from the root chosen internally by union--find,
so that union by rank remains available.

\begin{figure}[htbp]
\centering
\begin{subfigure}[b]{0.48\textwidth}
\centering
\includegraphics[width=0.8\linewidth]{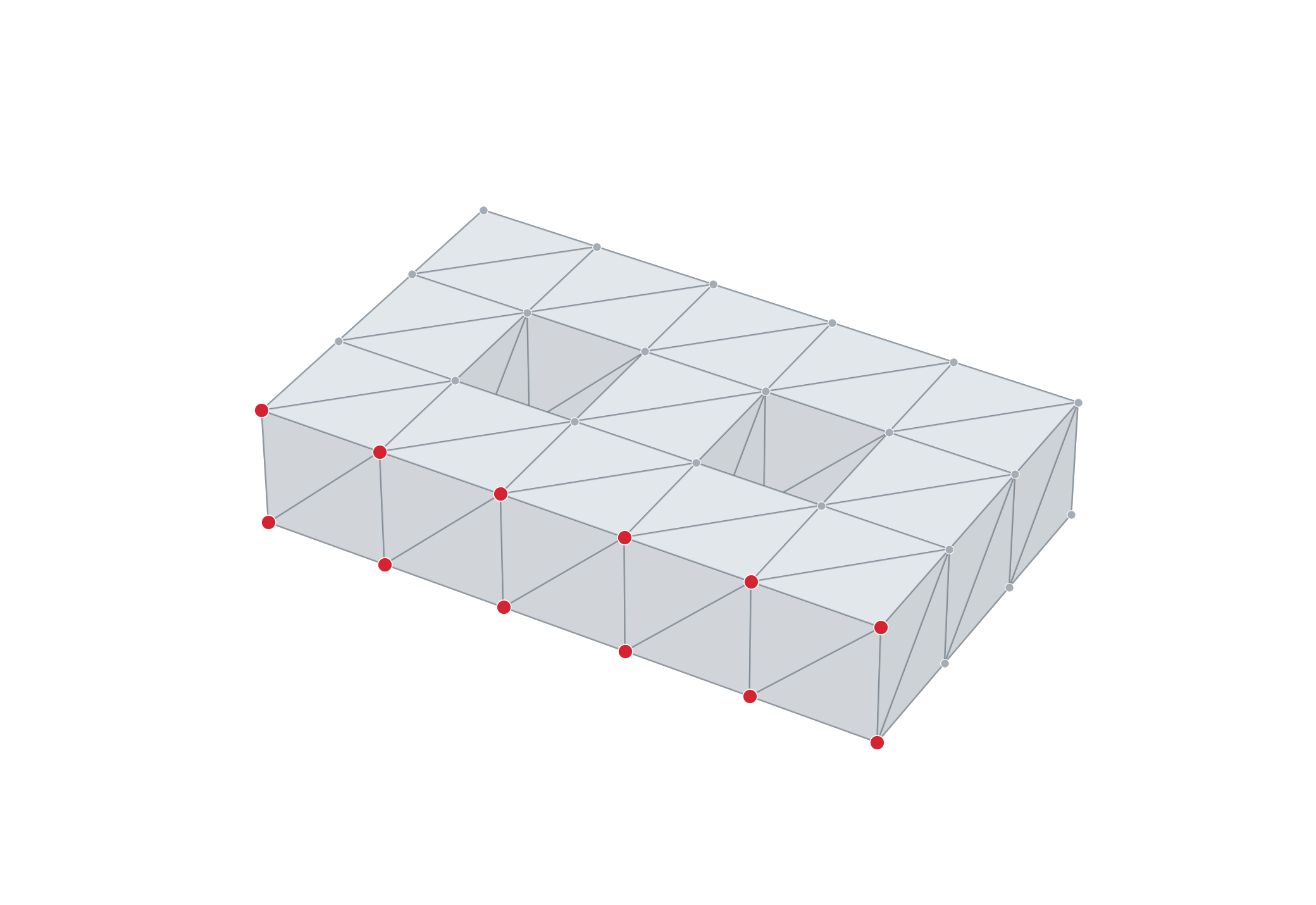}
\caption{Insert the vertices at the initial height.}
\label{fig:h01}
\end{subfigure}
\hfill
\begin{subfigure}[b]{0.48\textwidth}
\centering
\includegraphics[width=0.8\linewidth]{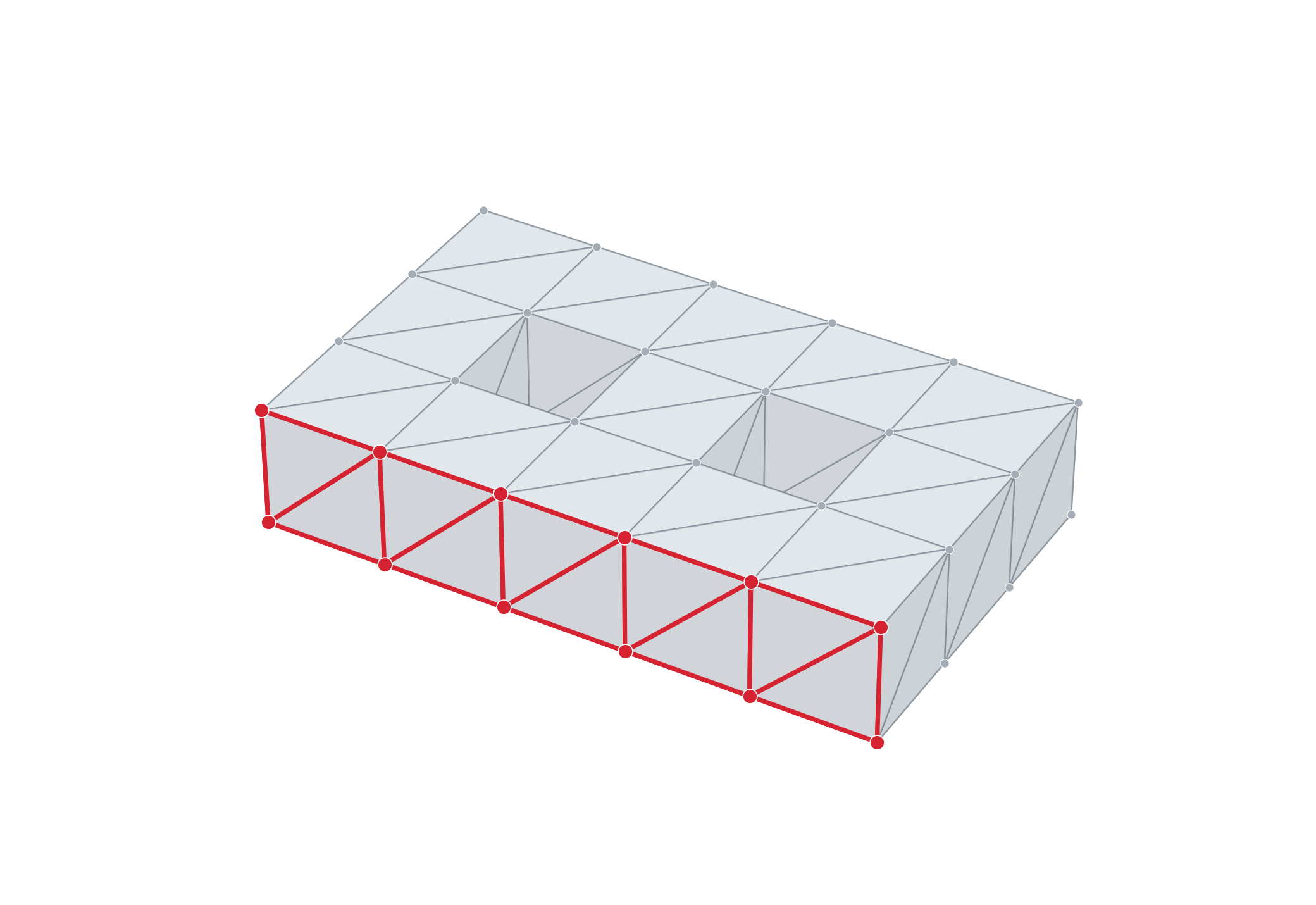}
\caption{Add the edges joining vertices at this height.}
\label{fig:h02}
\end{subfigure}

\vspace{0.3cm}

\begin{subfigure}[b]{0.48\textwidth}
    \centering
    \includegraphics[width=0.8\linewidth]{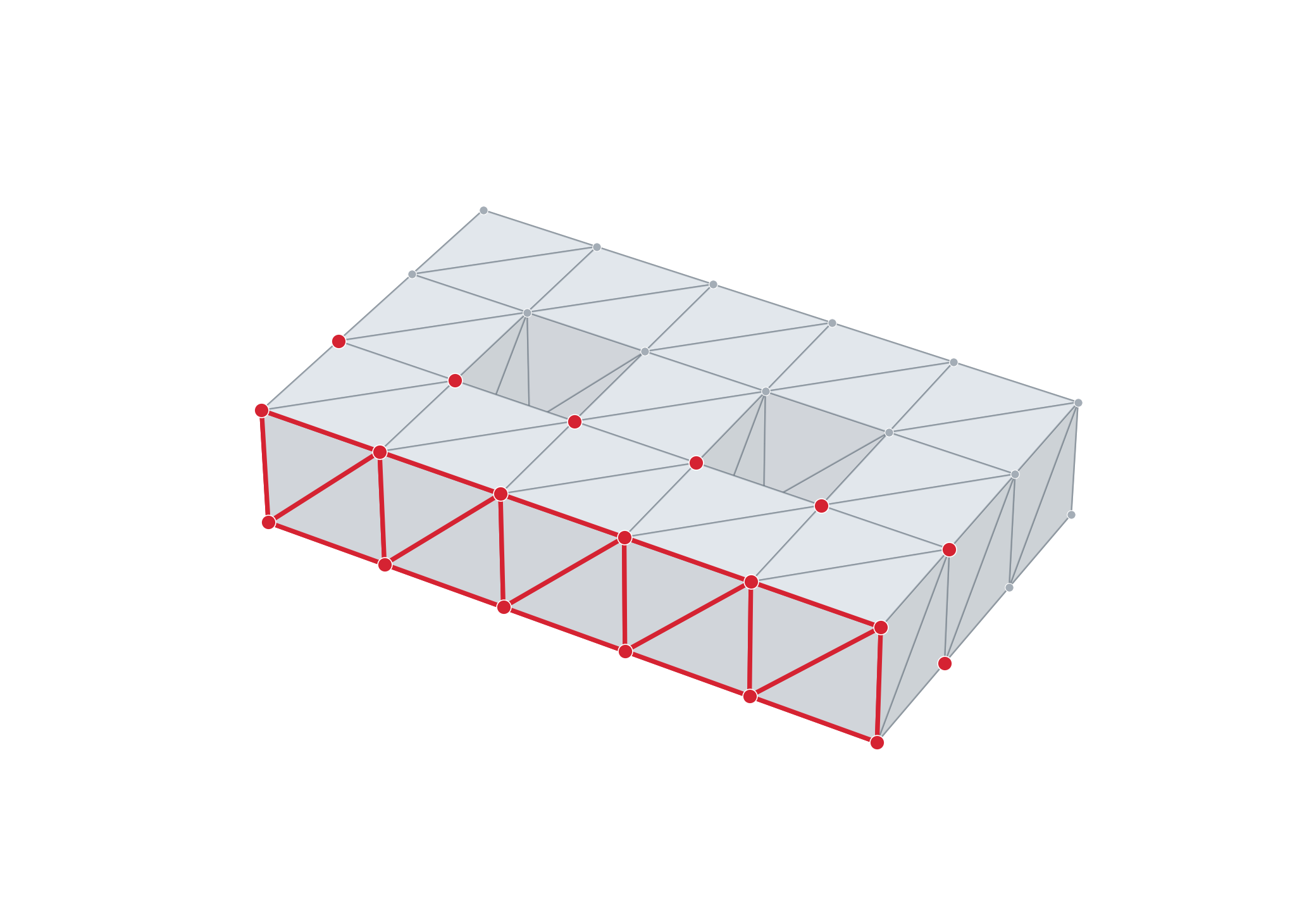}
    \caption{Insert the vertices at the next height.}
    \label{fig:h03}
\end{subfigure}
\hfill
\begin{subfigure}[b]{0.48\textwidth}
    \centering
    \includegraphics[width=0.8\linewidth]{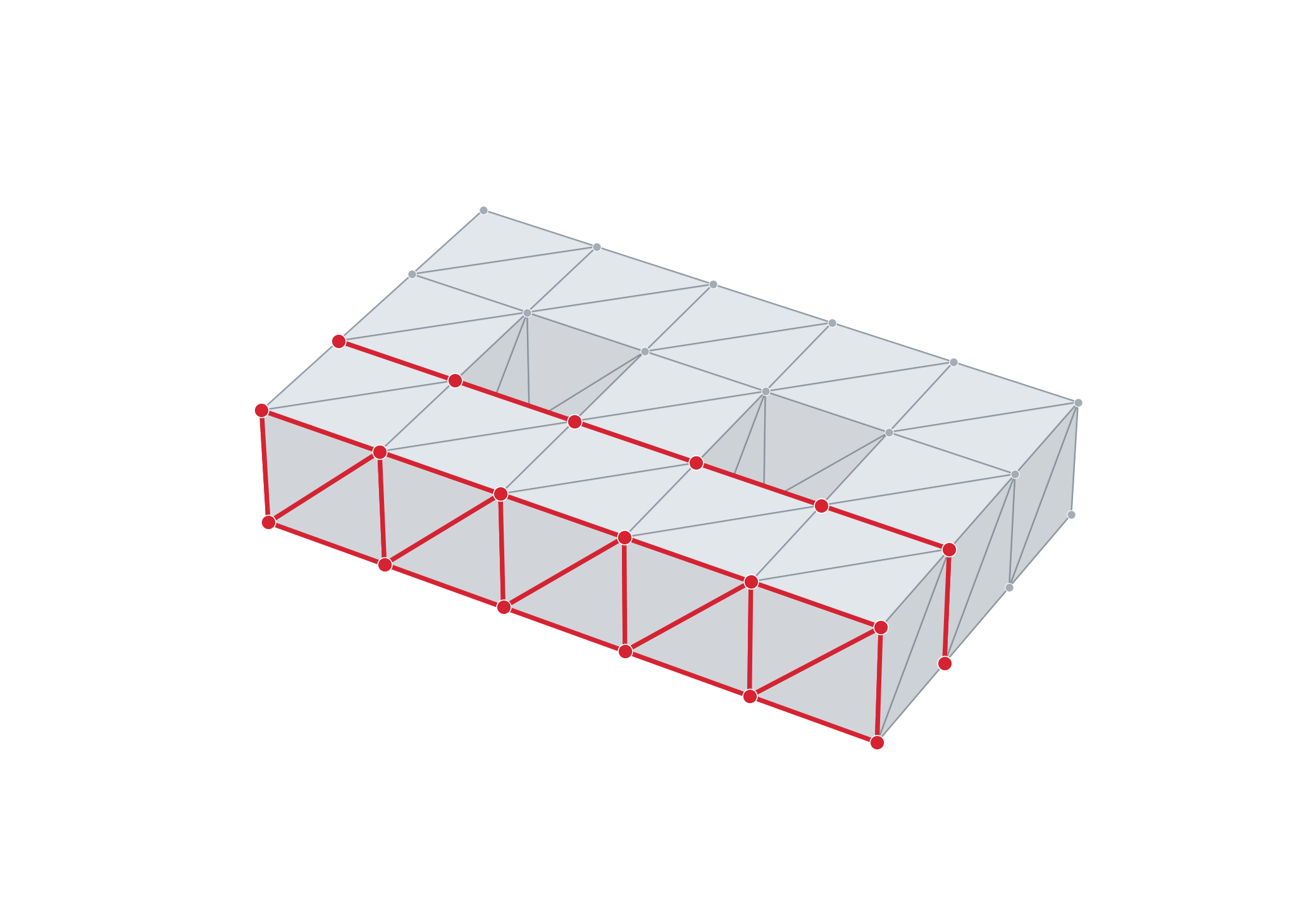}
    \caption{Add the edges joining the newly inserted vertices.}
    \label{fig:h04}
\end{subfigure}

\vspace{0.3cm}

\begin{subfigure}[b]{0.48\textwidth}
    \centering
    \includegraphics[width=0.8\linewidth]{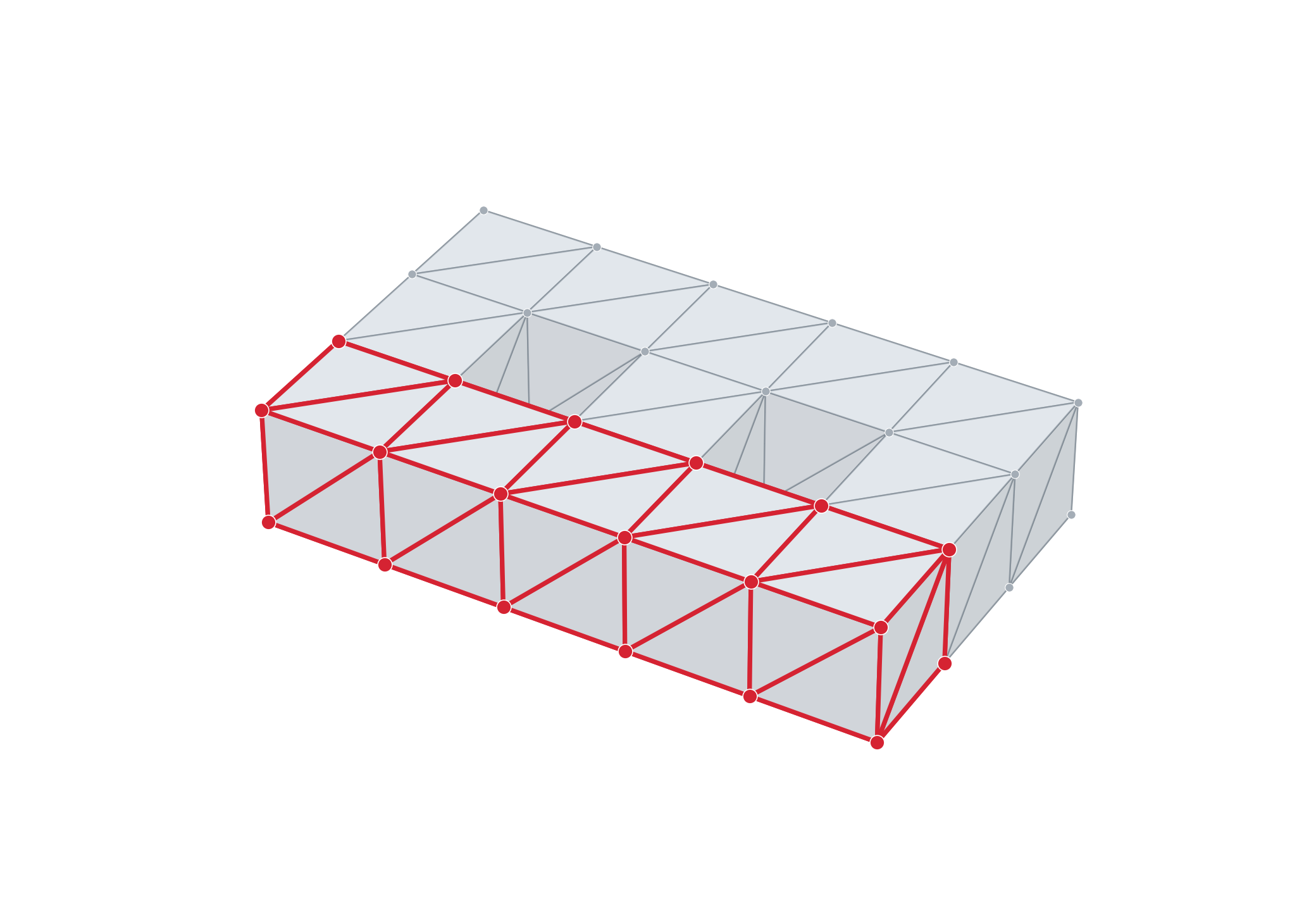}
    \caption{Add the edges connecting new vertices to lower vertices.}
    \label{fig:h05}
\end{subfigure}
\hfill
\begin{subfigure}[b]{0.48\textwidth}
    \centering
    \includegraphics[width=0.8\linewidth]{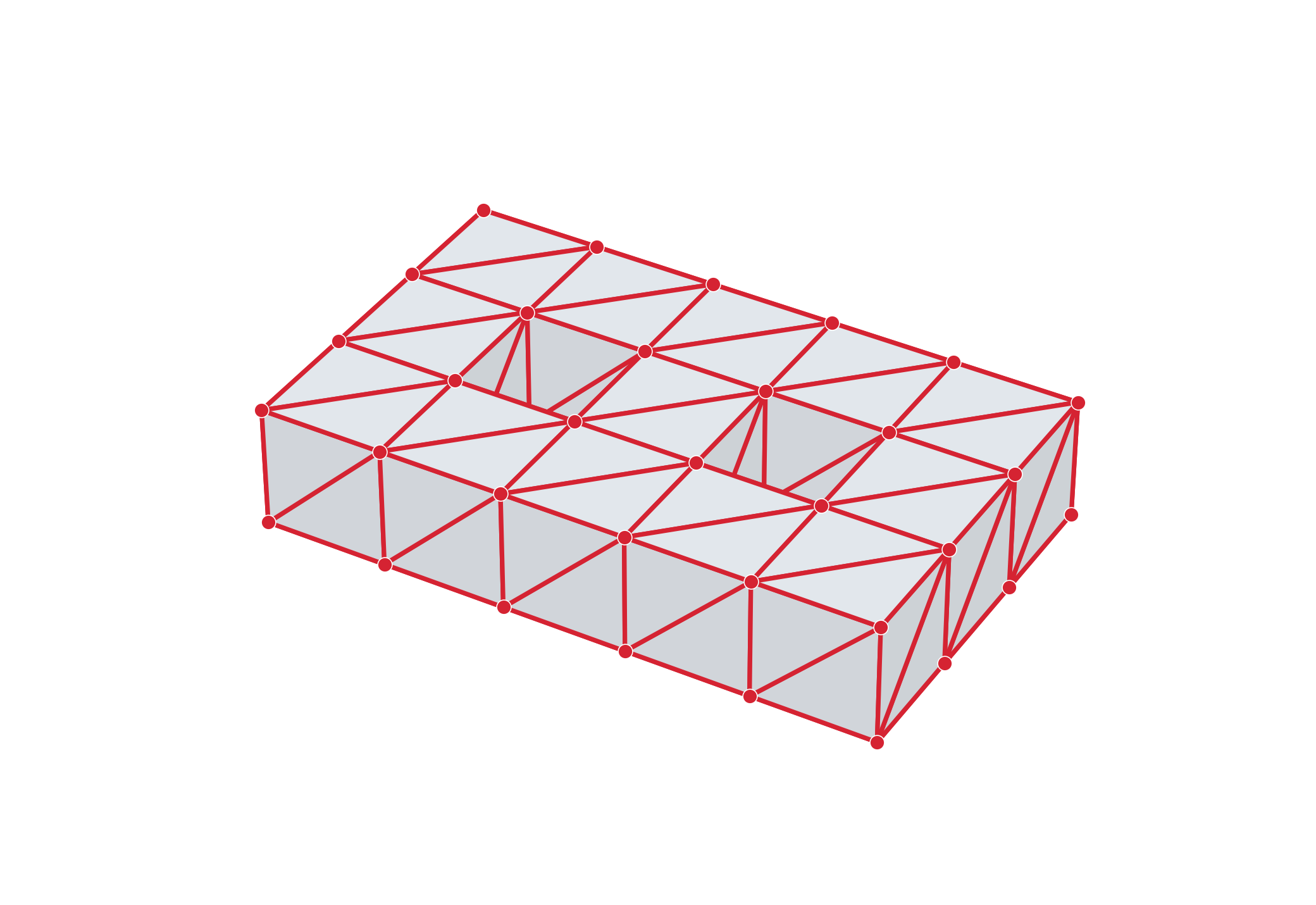}
    \caption{Continue through the filtration, recording component mergers.}
    \label{fig:h06}
\end{subfigure}

\caption{Computing degree-zero persistence by tracking connected components through the height filtration. At each merger, the younger component dies according to the elder rule.}
\label{fig:coho_svgs}

\end{figure}

Let $t_0$ be the smallest vertex height. We first insert all
vertices at height $t_0$
(Figure~\ref{fig:h01}), then add the edges whose endpoints
both have height $t_0$ (Figure~\ref{fig:h02}).

At each subsequent height $t$, the computation proceeds in
two stages. First, we insert the vertices at height $t$ and
the edges joining pairs of these vertices
(Figures~\ref{fig:h03} and~\ref{fig:h04}). Second, we add
the edges joining these new vertices to vertices at lower
heights (Figure~\ref{fig:h05}). Together, these stages insert
exactly the edges whose filtration value is $t$.

Whenever an edge joins two distinct components, we apply
the \textbf{elder rule}: the component with the earlier
birth survives, and the younger component dies at height
$t$. If both components have the same birth height, their
birth vertices are compared using the tie-breaking order.
A merger therefore produces an ordinary interval

$$
[b_{\mathrm{young}},t),
$$

provided that $b_{\mathrm{young}}<t$. An edge whose endpoints
already belong to the same component produces no change
in degree-zero homology.

We repeat this process until every vertex and edge has
been inserted (Figure~\ref{fig:h06}).

\subsection{Essential Intervals}

The components that remain after the final filtration
level correspond to the essential degree-zero classes.
For each connected component $C$ of $K$, its ordinary
persistence interval is

$$
\left[\min_{p\in V(C)}h_v(p),\infty\right).
$$

In extended persistence, this class dies when the descending
superlevel set first meets $C$, at its maximum height.
Consequently,

$$
\Ess_0(K,h_v)
=
\left\{
\left[
\bigl(\min_{p\in V(C)}h_v(p),\mathrm{ord}\bigr),
\bigl(\max_{p\in V(C)}h_v(p),\mathrm{rel}\bigr)
\right)
:
C\in\pi_0(K)
\right\},
$$

with intervals counted with multiplicity.

Thus the same connectivity computation determines both the ordinary intervals and the connected components needed to construct the essential intervals.

The following pseudocode returns the ordinary and essential intervals. The operation \textsc{MakeSet} inserts a singleton component into an initially empty union--find structure.
The array $\beta$ stores the birth vertex of each current component, independently of its union--find root.
Both output collections are multisets.

\begin{algorithm}[htbp]
\caption{\textsc{DegreeZeroExtendedPersistence}}
\begin{algorithmic}[1]
\Require Vertices $V$, edges $E$, height function $h:V\to\mathbb{R}$
\Ensure Ordinary intervals $\mathcal{O}$ and essential intervals $\mathcal{S}$

\State $U\gets$ an empty union--find structure
\State $\mathcal{O}\gets\varnothing$; $\mathcal{S}\gets\varnothing$
\State $L\gets$ the distinct vertex heights in increasing order

\For{each $a\in L$}
\For{each vertex $p$ with $h(p)=a$, in tie-breaking order}
\State $U.\Call{MakeSet}{p}$
\State $\beta[p]\gets p$
\EndFor

\State $E_a\gets\{[p,q]\in E:\max\{h(p),h(q)\}=a\}$
\State Order $E_a$ increasingly by $\lambda(e)$
\For{each $[p,q]\in E_a$}
    \State $r_p\gets U.\Call{Find}{p}$
    \State $r_q\gets U.\Call{Find}{q}$

    \If{$r_p\neq r_q$}
        \State $b_p\gets\beta[r_p]$; $b_q\gets\beta[r_q]$
        \State $(b_{\mathrm{old}},b_{\mathrm{young}})
        \gets$ $(b_p,b_q)$ sorted by $\kappa$

        \If{$h(b_{\mathrm{young}})<a$}
            \State Append $[h(b_{\mathrm{young}}),a)$ to $\mathcal{O}$
        \EndIf

        \State $U.\Call{Union}{r_p,r_q}$
        \State $r\gets U.\Call{Find}{p}$
        \State $\beta[r]\gets b_{\mathrm{old}}$
        \Comment{Preserve the older birth}
    \EndIf
\EndFor

\EndFor

\State $R\gets\{U.\Call{Find}{p}:p\in V\}$
\For{each $r\in R$}
\State $m[r]\gets-\infty$
\EndFor

\For{each $p\in V$}
\State $r\gets U.\Call{Find}{p}$
\State $m[r]\gets\max\{m[r],h(p)\}$
\EndFor

\For{each $r\in R$}
\State Append
$[(h(\beta[r]),\mathrm{ord}),(m[r],\mathrm{rel}))$
to $\mathcal{S}$
\EndFor

\State \Return $(\mathcal{O},\mathcal{S})$
\end{algorithmic}
\end{algorithm}

\section{Essential $H_1$ Classes}\label{chapter-h1classes}

We fix some notation here. Let $T$ be a tree contained in a larger graph. If $v_1,v_2\in T$ are the two endpoints of an edge $e\notin T$, then the \emph{unique} path in $T$ joining $v_1$ and $v_2$ is denoted by $P_T(e)$. 

\subsection{Tree--Cotree Structure}

We use a \textbf{tree--cotree decomposition} to construct cycle representatives whose homology classes form a basis of $H_1(K;\mathbb{F}_2)$ \cite{Eppstein2003,Erickson2012Cycles}. We first describe the construction for a connected surface; for a disconnected surface, it is applied separately to each connected component.

Fix a total order $\prec$ on $E$ as in \ref{total-order}. Construct a primal spanning tree $T$ by processing edges in increasing order and accepting an edge whenever it joins distinct components. Construct a dual spanning tree by processing the edges of $E\setminus T$ in the reverse order. Vertices in the dual tree are in bijection with the faces of the original graph. An edge in the dual tree is given by edges in the original graph, when two faces share a common edge. We restrict to edges in $E\setminus T$. Set
\[
C=\{e\in E:e^*\in T^*\},
\qquad
L=E\setminus(T\cup C).
\]

For each $e\in L$, let $\gamma_e$ be the cycle consisting of $e$ and $P_T(e)$. The tree--cotree basis theorem states that $\{[\gamma_e]:e\in L\}$ is a basis of $H_1(K;\mathbb{F}_2)$.

For a connected closed orientable surface of genus $g$,

$$
|L|
=
|E|-(|V|-1)-(|F|-1)
=
2g,
$$

where the final equality follows from the Euler characteristic formula.

\begin{figure}[htbp]
\centering
\begin{subfigure}[b]{0.48\textwidth}
    \centering
    \includegraphics[width=0.85\linewidth]{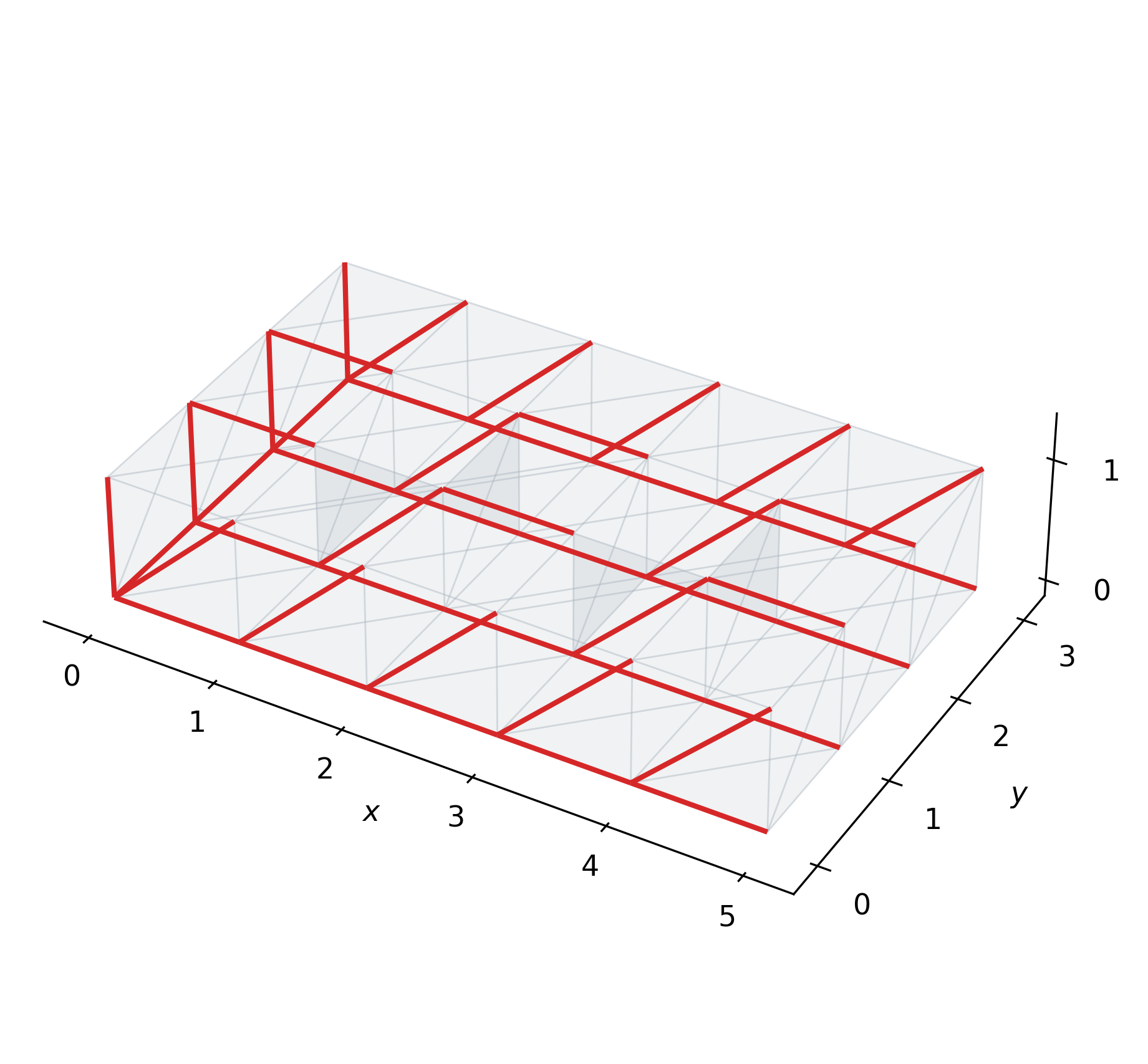}
    \caption{A triangulated genus-two surface with tree edges in red.}
    \label{fig:genus-two-tree}
\end{subfigure}
\hfill
\begin{subfigure}[b]{0.48\textwidth}
    \centering
    \includegraphics[width=0.85\linewidth]{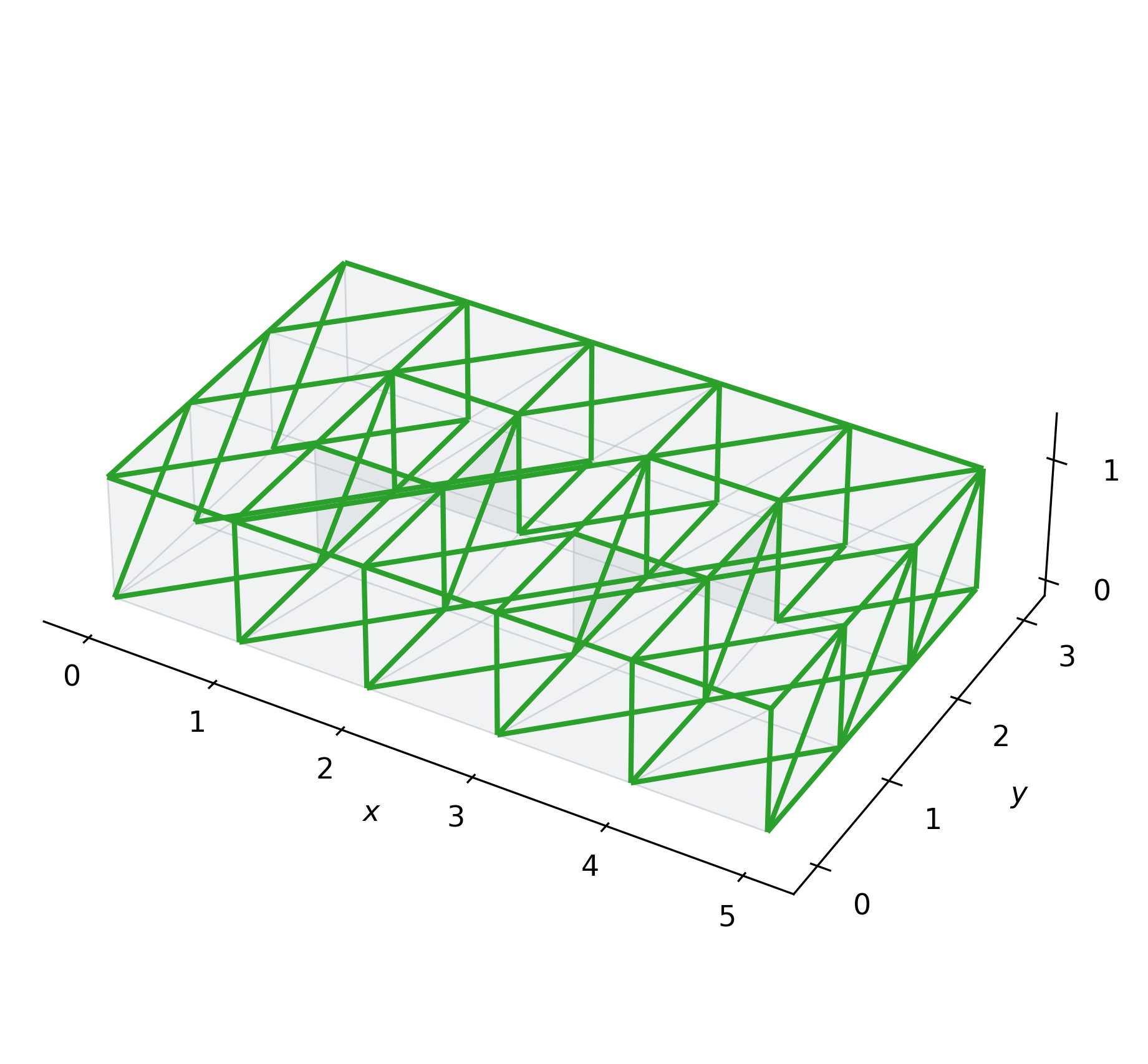}
    \caption{A triangulated genus-two surface with cotree edges in green.}
    \label{fig:genus-two-cotree}
\end{subfigure}

\vspace{0.3cm}

\begin{subfigure}[b]{0.48\textwidth}
    \centering
    \includegraphics[width=0.85\linewidth]{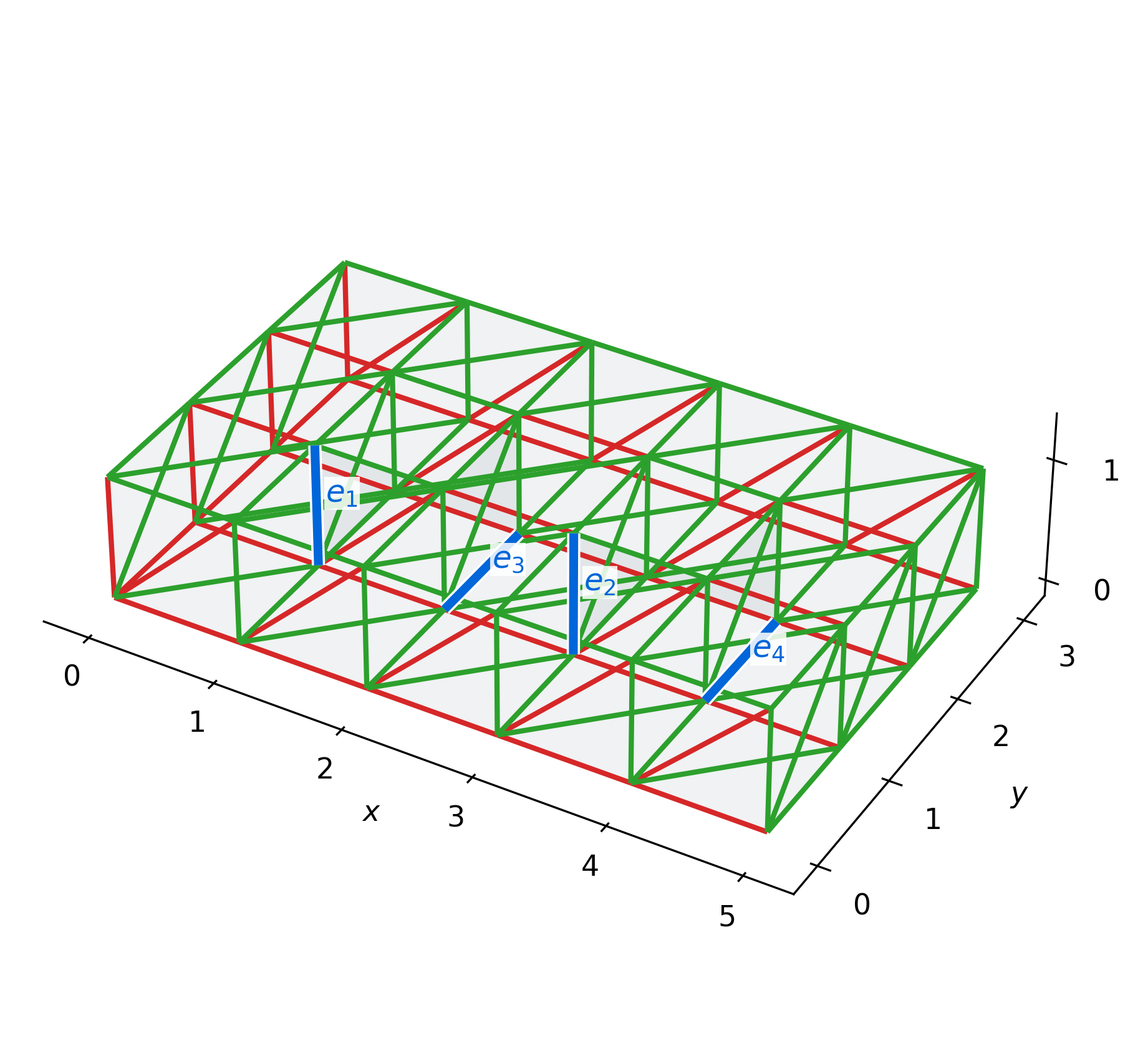}
    \caption{A triangulated genus-two surface with tree edges in red, cotree edges in green, and essential cycle-forming edges in blue.}
    \label{fig:genus-two-all}
\end{subfigure}

\caption{The three edge sets in the tree-cotree decomposition.}
\label{fig:tree-cotree}

\end{figure}

\begin{example}
    We triangulate a genus $2$ surface (Figure \ref{fig:genus-two-y}). Figure \ref{fig:genus-two-tree} shows the fundamental tree, while Figure \ref{fig:genus-two-cotree} shows the cotree structure. Note that the cotree structure is not a ture tree with respect to vertices, but it is indeed with respect to faces. Figure \ref{fig:genus-two-all} displays the leftover edges after combining the tree and the cotree structures. These leftover edges are the foundation of the essential $1$ cycles.
\end{example}

\begin{definition}
\label{def:image-filtration-compatible}
Let $K$ be a finite simplicial complex with vertex heights
$h\colon V(K)\to\mathbb{R}$. For each $t\in\mathbb{R}$, define
\[
    K_t
    =
    \{\sigma\in K : \max_{v\in\sigma} h(v)\le t\},
\]
and let
\[
    I_t
    =
    \operatorname{im}\!\left(
        H_1(K_t;\mathbb{F}_2)
        \longrightarrow
        H_1(K;\mathbb{F}_2)
    \right),
\]
where the map is induced by inclusion.

An ordered basis
$\mathcal{B}=([\gamma_1],\ldots,[\gamma_r])$
of $H_1(K;\mathbb{F}_2)$, together with assigned birth heights
$b_1\le\cdots\le b_r$, is \textbf{compatible with the image filtration}
if, for every $t\in\mathbb{R}$,
\[
    I_t
    =
    \operatorname{span}_{\mathbb{F}_2}
    \{[\gamma_i] : b_i\le t\}.
\]
Equivalently, $\{[\gamma_i]:b_i\le t\}$ is a basis of $I_t$
for every $t$. Equal birth heights are permitted.
\end{definition}

Assume we have obtained the tree-cotree decomposition $E=T\sqcup C\sqcup L$ as above. We show that it yields a basis compatible with the image filtration. Define
$$\gamma_e=e+P_T(e)$$
for every $e\in E\setminus T$, including the edges in $C$.

\begin{lemma}\label{fundamental-cycle-identity}
    With notation as above, 
\[
z=\sum_{e\in\operatorname{supp}(z)\setminus T}\gamma_e
\qquad
\text{for every }z\in Z_1(K).
\]
\end{lemma}

\begin{proof}
    Let $w=z+\sum_{e\in\operatorname{supp}(z)\setminus T}\gamma_e$. Since all summands are cycles, $\partial w=0$. Let $f\in \operatorname{supp}(w)$. If $f\in E\setminus T$, then $\gamma_f$ is the only cycle containing $f$ in the sense that if $f\in\gamma_i$, then $i=f$. There are two choices, $f\in \operatorname{supp}(z)\setminus T$ or not. For the first case, we have that $f\in z$ too, so that $f$ shows exactly twice on the right-hand side, and $f+f=0$, a contradiction to $f$ being in $\operatorname{supp}(w)$. For the latter case, if $f\notin \operatorname{supp}(z)$, then we have no copies of $f$ at all on the right-hand side of the equation, another contradiction. So $f\in T$. Therefore $\operatorname{supp}(w)\subseteq T$. Such a cycle must be zero \cite[Theorem~1.9.6]{Diestel2005GraphTheory}.
\end{proof}

\begin{theorem}\label{tree-cotree-theorem}
    Let $K=(V,E,F)$ triangulate a connected closed orientable surface, and work over $\mathbb F_2$. Given vertex heights $h:V\to\mathbb R$, define
\[
h_v(\sigma)=\max_{p\in\sigma}h(p),
\qquad
K_t=\{\sigma\in K:h_v(\sigma)\le t\}.
\]
Let $E=T\sqcup C\sqcup L$ be a tree-cotree decomposition obtained exactly as above. For each $\ell\in L$, let
\[
\gamma_\ell=\ell+P_T(\ell),
\qquad
b_\ell=h_v(\ell).
\]

Then the tree–cotree homology basis is compatible with the image filtration: for every $t\in\mathbb R$,
\[
\operatorname{im}\!\left(H_1(K_t)\longrightarrow H_1(K)\right)
=
\operatorname{span}
\{[\gamma_\ell]:\ell\in L,\ b_\ell\le t\}.
\]
In particular, $b_\ell$ is the birth height of $[\gamma_\ell]$ in this image filtration.
\end{theorem}

\begin{proof}
    The tree–cotree basis theorem gives that the classes
$\{[\gamma_\ell]:\ell\in L\}$ form a basis of $H_1(K)$. We prove the additional filtration property.

First consider the forward primal sweep. When an edge $e\notin T$ is processed, its endpoints are already connected by previously accepted tree edges. This path remains the unique path between its endpoints in the final tree $T$. Consequently, every edge of $P_T(e)$ precedes $e$, and therefore

\begin{equation}\label{eq:1}
    \gamma_e\in Z_1(K_{h_v(e)}).
\end{equation}

Now fix $c\in C$. Removing $c^*$ from from the dual tree partitions the faces into two nonempty sets $F_c^1$ and $F_c^2$. Without loss of generality, choose $F_c^1$ and define
\[
W_c=\sum_{\tau\in F_c^1}\tau,
\qquad
B_c=\partial W_c.
\]
Note that 
\begin{equation}\label{eq:2}
    \operatorname{supp}(B_c)\setminus T
=\{c\}\cup\Lambda_c
\end{equation}
for some subset $\Lambda_c\subseteq L$. Indeed, if $f\in \operatorname{supp}(B_c) \cap C$, then either $f$ is between two faces in $W_c$ or $f=c$. But $f\in\partial W_c$, so it follows that $f=c$. The only other possible edges in 
$ \operatorname{supp}(B_c)\setminus T$ are edges in $L$.

We claim that
\begin{equation}\label{eq:3}
    \ell\prec c
\qquad\text{for every }\ell\in\Lambda_c.
\end{equation}

Suppose instead that $c\prec\ell$. In the reverse dual sweep, $\ell^*$ would then be processed before $c^*$. At that moment, the accepted dual forest is a subgraph of $T^*\setminus\{c^*\}$ which we know to be partitioned in two. Because $\ell^*$ joins the two components of this latter forest, its endpoints cannot yet be connected. The algorithm would therefore accept $\ell^*$, contradicting $\ell\in L$. This proves \eqref{eq:3}.

Next, applying Lemma \ref{fundamental-cycle-identity} and \eqref{eq:2} to $B_c$ directly yields
\[
B_c=\gamma_c+\sum_{\ell\in\Lambda_c}\gamma_\ell.
\]
Since $B_c=\partial W_c$ is a boundary in $K$, we obtain

\begin{equation}\label{eq:4}
    [\gamma_c]
=
\sum_{\ell\in\Lambda_c}[\gamma_\ell]
\quad\text{in }H_1(K),
\qquad
h_v(\ell)\le h_v(c).
\end{equation}

Thus every fundamental cycle associated with a cotree edge is homologous in $K$ to a sum of fundamental cycles associated with leftover edges, each appearing no later in the filtration than the original cycle.

We can now prove the claimed equality. By \eqref{eq:1}, every $\gamma_\ell$ with $b_\ell\le t$ is a cycle in $K_t$, so
\[
\operatorname{span}\{[\gamma_\ell]:b_\ell\le t\}
\subseteq
\operatorname{im}\!\left(H_1(K_t)\to H_1(K)\right).
\]
Conversely, let $z\in Z_1(K_t)$. Using Lemma \ref{fundamental-cycle-identity}, we can express $z$ using fundamental cycles $\gamma_e$ whose closing edges satisfy $h_v(e)\le t$. For each such $e\in C$, apply \eqref{eq:4}. That is:
\[
b_\ell=h_v(\ell)\le h_v(e)\le t.
\]
Hence
\[
[z]\in\operatorname{span}\{[\gamma_\ell]:b_\ell\le t\},
\]
which proves the reverse inclusion.

Finally, $[\gamma_\ell]$ belongs to the image at $t=b_\ell$. For $t<b_\ell$, linear independence prevents the image from containing $[\gamma_\ell]$ (recall that the tree-cotree theorem gives $[\gamma_\ell]_{\ell\in L}$ as a basis). Its birth height is therefore exactly $b_\ell$. 
\end{proof}

The following pseudo-code constructs a homology basis using the ordered tree--cotree decomposition. The candidate edges are stored in a list so that their processing order is preserved. For each returned cycle, the height of its closing edge is the height at which that particular representative appears. 

\begin{algorithm}[htbp]
\caption{\textsc{TreeCotreeEssentialCycles}}
\begin{algorithmic}[1]
\Require A closed surface mesh $K=(V,E,F)$ and vertex heights $h$
\Ensure Cycle representatives forming a basis of $H_1(K)$

\State $T\gets\varnothing$
\Comment{Primal spanning forest}
\State $\mathcal{C}\gets$ empty ordered list
\State $U_T\gets\Call{UnionFind}{V}$
\State Order edges increasingly by $\lambda(e)$

\For{each edge $e=[p,q]$ in this order}
\If{$U_T.\Call{Find}{p}=U_T.\Call{Find}{q}$}
\State Append $e$ to $\mathcal{C}$
\Else
\State $T\gets T\cup\{e\}$
\State $U_T.\Call{Union}{p,q}$
\EndIf
\EndFor

\State $U^*\gets\Call{UnionFind}{F}$
\State $L\gets$ empty list

\For{each $e\in\mathcal{C}$ in reverse order}
\State Let $f_1,f_2$ be the two faces incident to $e$
\If{$U^*.\Call{Find}{f_1}=U^*.\Call{Find}{f_2}$}
\State Append $e$ to $L$
\Else
\State $U^*.\Call{Union}{f_1,f_2}$
\EndIf
\EndFor

\State Restore the original edge order on $L$
\For{each $e=[p,q]\in L$}
\State $P_e\gets$ the unique path from $p$ to $q$ in $T$
\State $\gamma_e\gets P_e+e$
\Comment{Sum of edge chains over $\mathbb{F}_2$}
\EndFor

\State \Return $(\gamma_e)_{e\in L}$
\end{algorithmic}
\end{algorithm}

\section{Essential $H_1$ Pairing}\label{chapter-pairing}

Throughout this subsection, $K$ triangulates the closed
boundary surface $A=\partial X$, and all computations are
performed over $\mathbb{F}_2$. 

\subsection{Constructing dual cocycles.}

Let $[\gamma_1],\ldots,[\gamma_r]$ be a basis of $H_1(K)$, with birth heights $b_1\le\cdots\le b_r$ such that it is compatible with the image filtration. For a connected surface of genus $g$, $r=2g$.

When birth heights coincide, we refer to the tie-breaking convention described above in Section \ref{total-order}. Note that all of this is guaranteed by Theorem \ref{tree-cotree-theorem}.

By Corollary~\ref{evaluation-dual-basis}, there are
cohomology classes
$[\alpha_1],\ldots,[\alpha_r]$ satisfying

$$
\alpha_i(\gamma_j)=\delta_{ij}.
$$

We construct cocycle representatives of these classes
by solving linear systems.

Let $E=\{e_1,\ldots,e_m\}$ be the edge set, and let
$B_2$ be the edge--face boundary matrix. By proposition~\ref{boundary_matrix}, a cochain with
coordinate vector $x\in\mathbb{F}_2^m$ is a cocycle
precisely when

$$
B_2^\top x=0.
$$

Equivalently, the cocycle $x$ evaluates to zero on the boundary
of every triangular face.

Let $G\in\mathbb{F}_2^{m\times r}$ be the matrix whose
$j$-th column is the edge-coordinate vector of
$\gamma_j$. Thus

$$
G_{kj}
=
\text{the coefficient of }e_k\text{ in }\gamma_j.
$$

The evaluation-dual cocycle $\alpha_i$ is represented
by any solution of

$$
\begin{bmatrix}
B_2^\top\\
G^\top
\end{bmatrix}
x^{(i)}
=
\begin{bmatrix}
0\\
\mathbf{e}_i
\end{bmatrix},
$$

where $\mathbf{e}_i$ is the $i$-th standard basis vector
of $\mathbb{F}_2^r$.

These systems are consistent because the classes
$[\gamma_j]$ form a homology basis. Their solutions
need not be unique, but different solutions for the
same right-hand side differ by a coboundary and therefore
represent the same cohomology class.

Writing

$$
C=
\begin{bmatrix}
x^{(1)}&\cdots&x^{(r)}
\end{bmatrix},
$$

we may solve for all cocycles simultaneously:

$$
\begin{bmatrix}
B_2^\top\\
G^\top
\end{bmatrix}C
=
\begin{bmatrix}
0\\
I_r
\end{bmatrix}.
$$

The resulting representatives are

$$
\alpha_i=\sum_{k=1}^m C_{ki}e_k^*.
$$
Here $e_k^*e_j=\delta_{kj}$
\subsection{Recovering the essential pairing.}

For a threshold $t$, define the upper-star subcomplex

$$
K^t=
\{\sigma\in K:h(p)\ge t
\text{ for every vertex }p\in\sigma\}.
$$

Its vertices and edges are denoted by $V^t$ and $E^t$.
As $t$ decreases, these subcomplexes grow. Let $B_1^t$ be the vertex--edge boundary matrix of $K^t$.
Let $E(K^t)=\{e_1,\ldots,e_{m_t}\}$ be the
edges of $K^t$. Define
\[
C_t=
\begin{pmatrix}
\alpha_1(e_1) & \alpha_2(e_1) & \cdots & \alpha_r(e_1)\\
\alpha_1(e_2) & \alpha_2(e_2) & \cdots & \alpha_r(e_2)\\
\vdots & \vdots & \ddots & \vdots\\
\alpha_1(e_{m_t}) & \alpha_2(e_{m_t}) & \cdots &
\alpha_r(e_{m_t})
\end{pmatrix}.
\]
Thus, the $i$-th column of $C_t$ contains the edge
values of the restricted cocycle $\alpha_i|_{K^t}$.

We give an algorithm to obtain the essential pairing via the evaluation-dual cocycles, and show its correctness.

\begin{proposition}\label{prop:evaluation-dual-pairing}
Let $[\gamma_1],\ldots,[\gamma_r]$ be a basis,
with birth heights $b_1\le\cdots\le b_r$, compatible with the image
filtration of Definition~\ref{def:image-filtration-compatible}.
Set $F_i=\operatorname{span}\{[\gamma_1],\ldots,[\gamma_i]\}$ and
$F_0=0$. Construct cocycles satisfying
$\alpha_i(\gamma_k)=\delta_{ik}$ as described above.

For each superlevel subcomplex $K^t$,
compute a basis for $\ker B_1^t$, and assemble them as columns for $Z_t$ (so $Z_t$ may be interpreted as a basis for the space $Z_1(K^t)$), let $C_t$ contain the edge values of
the restricted cocycles, and set
\[
Q_t=C_t^{\mathsf T}Z_t,
\qquad
\rho_i(t)=i-\operatorname{rank}Q_t+\operatorname{rank}Q_t^{>i},
\]
where $Q_t^{>i}$ consists of rows $i+1,\ldots,r$.
Then, for the natural map $q_t:H_1(K)\longrightarrow H_1(K,K^t)$,
\begin{equation}\label{eq:evaluation-dual-relative-rank}
\rho_i(t)=\operatorname{rank}(q_t|_{F_i}).
\end{equation}

Let $t_1>\cdots>t_N$ be the distinct vertex heights, evaluating
each relative stage after inserting all simplices entering at
that height. Write $\rho_{i,j}=\rho_i(t_j)$, set
$\rho_{i,0}=i$, and define
\[
\Delta_{i,j}=\rho_{i,j-1}-\rho_{i,j},
\qquad
\mu_{i,j}=\Delta_{i,j}-\Delta_{i-1,j}.
\]
The intervals
\[
\bigl[(b_i,\mathrm{ord}),(t_j,\mathrm{rel})\bigr)
\qquad\text{for all }(i,j)\text{ with }\mu_{i,j}=1
\]
are exactly the essential degree intervals. In particular, $\mu_{i,j}\in\{0,1\}$.
\end{proposition}

\begin{proof}
Consider the evaluation mapping
\begin{align*}
    \Theta:H_1(K) &\to\mathbb F_2^r\\
    [z] &\mapsto\bigl(\alpha_1(z),\ldots,\alpha_r(z)\bigr)^{\mathsf T}.
\end{align*}
It is not hard to show it is a well-defined isomorphism of vector spaces.

Set $S_t=\operatorname{im}(H_1(K^t)\longrightarrow H_1(K))$.
Since the columns of $Z_t$ span all cycles of $K^t$, the columns
of $Q_t$ span their global homology coordinates. Thus
\begin{equation}\label{eq:evaluation-dual-superlevel-image}
\operatorname{col}Q_t=\Theta(S_t).
\end{equation}
Exactness of the long exact sequence of the pair $(K,K^t)$ gives
$\ker q_t=S_t$. $\ker (q_t|_{F_i})=F_i\cap\ker(q_t)=F_i\cap S_t$, therefore by rank-nullity
\[
\operatorname{rank}(q_t|_{F_i})
=i-\dim(F_i\cap S_t).
\]
The projection discarding the first $i$ coordinates has kernel
$\Theta(F_i)$. Its restriction to $\operatorname{col}Q_t$ has
rank $\operatorname{rank}Q_t^{>i}$. A similar argument as above gives
\[
\dim(F_i\cap S_t)
=\operatorname{rank}Q_t-\operatorname{rank}Q_t^{>i},
\]
and combining this with the previous expression proves \eqref{eq:evaluation-dual-relative-rank}.

We now show directly that the births and deaths agree.
Write $S_j=S_{t_j}$ and $S_0=0$. These spaces satisfy
\[
0=S_0\subseteq S_1\subseteq\cdots\subseteq S_N=H_1(K).
\]
For $i\ge1$, define
\[
\delta_{i,j}=\rho_{i,j}-\rho_{i-1,j}
=\dim\frac{F_i+S_j}{F_{i-1}+S_j}.
\]
Because $F_i/F_{i-1}$ is one-dimensional,
$\delta_{i,j}\in\{0,1\}$. Moreover, $\delta_{i,j}=0$ is
equivalent to
\[
[\gamma_i]\in F_{i-1}+S_j.
\]
Once this holds, it continues to hold as $j$ increases.
Since $\delta_{i,0}=1$ and $\delta_{i,N}=0$, there is a unique
step $j_i$ for which $\delta_{i,j_i -1} = 1$ and $\delta_{i,j_i} = 0$. 

By the definition of $j_i$, choose $a_{ik}$ so that 
\[
\eta_i:=[\gamma_i]+\sum_{k<i}a_{ik}[\gamma_k]\in S_{j_i}.
\]
This triangular change of basis preserves births in the following sense:
\[
F_i=\operatorname{span}\{\eta_1,\ldots,\eta_i\}.
\]
In particular, $\eta_i$ has birth height $b_i$. Moreover, by construction:
\[
S_j=\operatorname{span}\{\eta_i:j_i\le j\}.
\]
Thus the same basis is compatible with both the ordinary birth
filtration and the relative kernel filtration.

To relate this to extended persistence, we note that by the first isomorphism theorem and naturality of the maps of pairs, the following diagram commutes

\[\begin{tikzcd}
	\cdots & {H_1(K)} & {H_1(K)/S_{1}} & {H_1(K)/S_{2}} & \cdots\\
    \cdots & {H_1(K)} & {H_1(K,K^{t_1})} & {H_1(K,K^{t_{2}})} & \cdots
	\arrow[from=1-1, to=1-2]
	\arrow[from=1-2, to=1-3]
	\arrow[from=1-2, to=2-2]
	\arrow[from=1-3, to=1-4]
	\arrow[from=1-3, to=2-3]
	\arrow[from=1-4, to=1-5]
	\arrow[from=1-4, to=2-4]
	\arrow[from=2-1, to=2-2]
	\arrow[from=2-2, to=2-3]
	\arrow[from=2-3, to=2-4]
	\arrow[from=2-4, to=2-5]
\end{tikzcd}\]

with the leftmost side of the diagram being identities. Therefore, if $\eta_i$ is born at $b_i$ we can track where it disappears in the upper sequence of the diagram, which happens to be at $t_{j_i}$ which is when it gets quotiented out. Therefore that corresponds to the death of $\eta_i$, and we obtain the interval [($b_i,\ord$),($t_{j_i},\rel))$.

Moreover, since at precisely $j_i$ $\delta_{i,-}$ becomes $0$,
\[
\mu_{i,j}
=\delta_{i,j-1}-\delta_{i,j}
=
\begin{cases}
1,&j=j_i,\\
0,&j\ne j_i.
\end{cases}
\]

\end{proof}

\begin{algorithm}[htbp]
\caption{\textsc{EssentialPairingFromDualCocycles}}
\begin{algorithmic}[1]
\Require Mesh $K=(V,E,F)$, heights $h$, a filtration-compatible homology basis with births $b_1,\ldots,b_r$, and its dual-cocycle matrix $C$
\Ensure The multiset $\mathcal{I}$ of essential degree-one intervals

\State $\mathcal{T}\gets$ distinct vertex heights in decreasing order
\State $\mathcal{I}\gets\varnothing$
\For{$i=0,\ldots,r$}
\State $\rho_i^{\mathrm{old}}\gets i$
\Comment{Initially the superlevel set is empty}
\EndFor

\For{each $t\in\mathcal{T}$}
\State $V^t\gets\{p\in V:h(p)\ge t\}$
\State $E^t\gets\{[p,q]\in E:p,q\in V^t\}$
\State Form the vertex--edge boundary matrix $B_1^t$
\State Compute a cycle basis $Z_t$ for $\ker B_1^t$
\State $C_t\gets C$ restricted to rows indexed by $E^t$
\State $Q_t\gets C_t^\top Z_t$
\State $\rho_0^{\mathrm{new}}\gets0$
\State $\Delta_0\gets0$

\For{$i=1,\ldots,r$}
    \State $\rho_i^{\mathrm{new}}
    \gets i-\operatorname{rank}(Q_t)
    +\operatorname{rank}(Q_t^{>i})$
    \State $\Delta_i\gets
    \rho_i^{\mathrm{old}}-\rho_i^{\mathrm{new}}$
    \State $\mu_i\gets\Delta_i-\Delta_{i-1}$
    \If{$\mu_i=1$}
        \State Append $[(b_i,\mathrm{ord}),(t,\mathrm{rel}))$
        to $\mathcal{I}$
    \EndIf
\EndFor

\State $\rho^{\mathrm{old}}\gets\rho^{\mathrm{new}}$

\EndFor

\State \Return $\mathcal{I}$
\end{algorithmic}
\end{algorithm}

The rank of a matrix with no rows or no columns is
taken to be zero. At the end of the sweep, $K^t=K$,
so $S_t=H_1(K)$ and every essential interval has died.

\begin{figure}[htbp]
    \centering
    \includegraphics[width=0.85\linewidth]{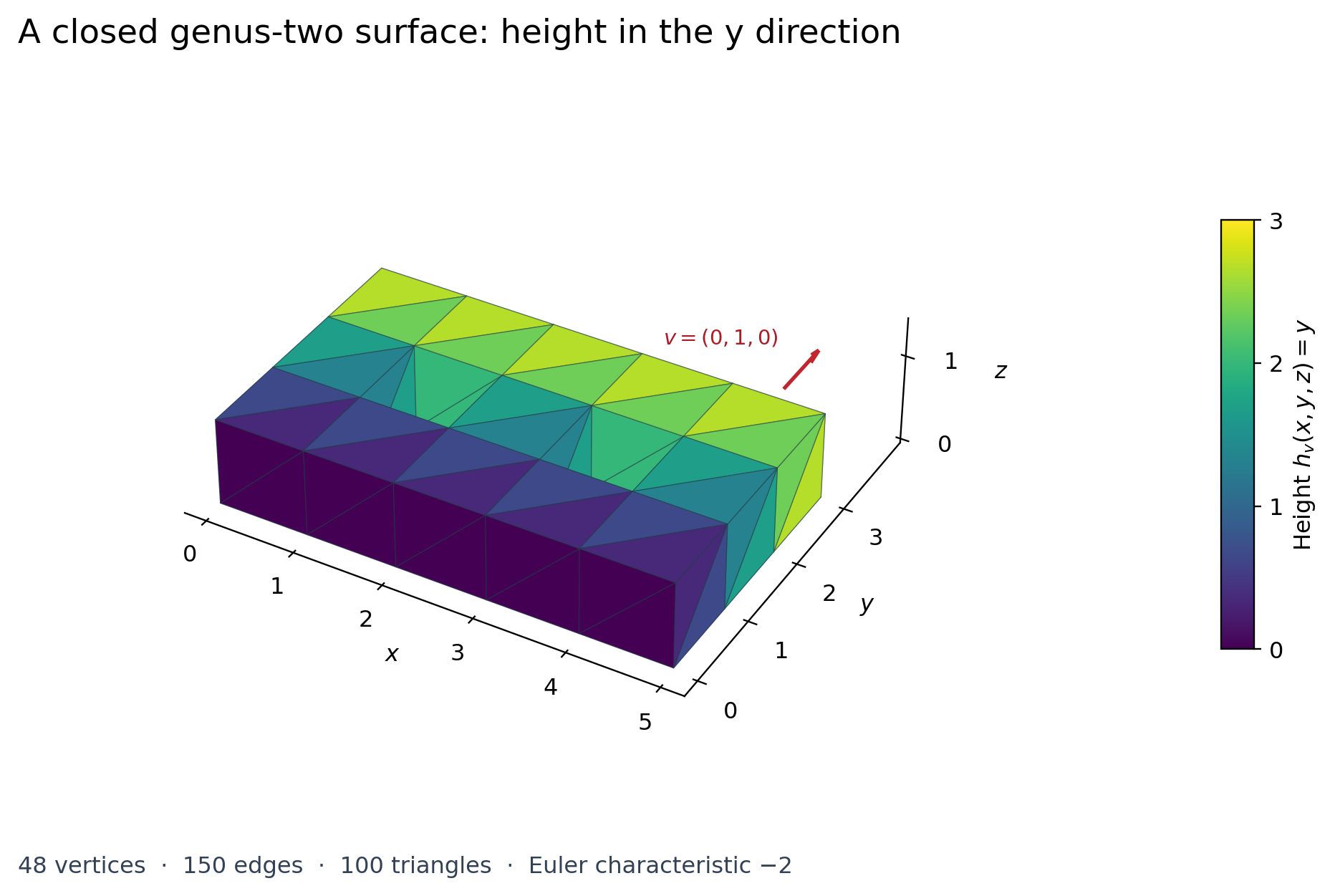}
    \caption{A triangulated genus-two surface coloured by
    the height function $h_v(x,y,z)=y$, for $v=(0,1,0)$.}
    \label{fig:genus-two-y}
\end{figure}

\begin{example}
    We will do a genus $2$ surface. Consider Figure \ref{fig:genus-two-y}. Let
\[
P=([0,5]\times[0,3])\setminus
\bigl(((1,2)\times(1,2))\cup((3,4)\times(1,2))\bigr),
\qquad X=P\times[0,1],\qquad A=\partial X.
\]
Triangulate the boundary as in the figure. There are
$48$ vertices, $150$ edges, and $100$ triangles, so
$\chi(A)=-2$ and $\dim H_1(A;\mathbb F_2)=4$.
Choose
\[
v=(0,1,0),\qquad h_v(x,y,z)=y.
\]
The vertex heights are $0,1,2,3$. At each numerical height,
insert all simplices at that height before computing ranks.
All calculations below are over $\mathbb F_2$.

Order the edges by $\lambda$ as in Section \ref{total-order}. The
ordered tree--cotree construction gives four fundamental
cycles with the following closing edges and birth heights:
\[
\begin{array}{ccc}
\toprule
\text{Cycle}&\text{Closing edge}&b_i\\
\midrule
\gamma_1&[(1,1,0),(1,1,1)]&1\\
\gamma_2&[(3,1,0),(3,1,1)]&1\\
\gamma_3&[(2,2,0),(2,1,0)]&2\\
\gamma_4&[(4,2,0),(4,1,0)]&2\\
\bottomrule
\end{array}
\]
Writing $H=H_1(A)$ and
$F_i=\operatorname{span}\{[\gamma_1],\ldots,[\gamma_i]\}$,
the ordinary image filtration is
\[
\operatorname{im}\bigl(H_1(A_a)\longrightarrow H\bigr)
=\begin{cases}
0,&a<1,\\
F_2,&1\le a<2,\\
H,&a\ge2.
\end{cases}
\]
Thus $(b_1,b_2,b_3,b_4)=(1,1,2,2)$.

Let $G$ be the matrix of these cycles, let $B_2$ be the
edge--face boundary matrix, and choose evaluation-dual
cocycles by solving
\[
\begin{bmatrix}B_2^{\mathsf T}\\G^{\mathsf T}\end{bmatrix}C
=\begin{bmatrix}0\\I_4\end{bmatrix}.
\]
At each superlevel stage, build $Z_t$  using linear algebra (compute a basis of $\ker (B_1^T)$).

The computation gives
\[
S_t=\begin{cases}
0,&t>2,\\
F_2,&1<t\le2,\\
H,&t\le1.
\end{cases}
\]
In particular, the opposite tunnel walls entering at
$y=2$ supply representatives of the two classes born at
$y=1$. Then, computing $Q_t=C_t^TZ_t$
\[
Q_2=\begin{pmatrix}
1&0\\
0&1\\
0&0\\
0&0
\end{pmatrix},\qquad Q_1=I_4.
\]

Using
\[
\rho_i(t)=i-\operatorname{rank}Q_t
+\operatorname{rank}Q_t^{>i},\qquad
\mu_i(t)=\Delta_i(t)-\Delta_{i-1}(t),
\]
where $\Delta_i(t)$ is the drop in $\rho_i$ at height $t$ and $Q_t^{>i}$ consists of rows $i+1,\ldots,r$.
gives
\[
\begin{array}{cccc}
\toprule
\text{Relative stage}&(\rho_1,\rho_2,\rho_3,\rho_4)
 &(\Delta_1,\Delta_2,\Delta_3,\Delta_4)
 &(\mu_1,\mu_2,\mu_3,\mu_4)\\
\midrule
t>2&(1,2,3,4)&---&---\\
t=2&(0,0,1,2)&(1,2,2,2)&(1,1,0,0)\\
t=1&(0,0,0,0)&(0,0,1,2)&(0,0,1,1)\\
\bottomrule
\end{array}
\]
There is no further essential rank change at $t=0$.
Consequently, the essential degree-one interval multiset
of the boundary consists of two copies of
\[
\bigl[(1,\mathrm{ord}),(2,\mathrm{rel})\bigr)
\]
and two copies of
\[
\bigl[(2,\mathrm{ord}),(1,\mathrm{rel})\bigr).
\]

For the solid, consider the left-over edges from the tree-cotree algorithm. $e_1$ and $e_2$ have normals $\frac{{(0,1,0)+(1,0,0)}}{|(0,1,0)+(1,0,0)|}$. This gives a negative sign with respect to the direction of the height function. On the other hand, for $e_3$ and $e_4$, the corresponding sum is $(0,-1,0)+(1,0,0)$. This yields a positive sign. Therefore these are the cycles being kept, and the essential
degree-one interval multiset of $X$ consists of two copies
of
\[
\bigl[(2,\mathrm{ord}),(1,\mathrm{rel})\bigr).
\]
\end{example}

\section{Conclusions and Future Work}

To summarize our work, we implemented an algorithm that computes a discrete version of the $3$D $\XPHT$. This involved reducing the problem to computing $\XPHT_0$ and $\Ess_1$ on a bounding surface. Computing $\Ess_1$ involved a technical step when pairing the critical points that provide the births and deaths of the essential classes. 

Experimentally, it is known that computing zeroth degree persistence often provides a ``good'' and efficient \cite{Turner2014PHT} summary of the object of study. Table \ref{tab:xpht} provides mathematical justification in the case of the $\XPHT$: one obtains most of the higher dimensional diagrams from just $\Ord_0$.

 In \cite{LewinerLopesTavares2003} and \cite{Bauer2011}, albeit with a different name, they use the tree-cotree decomposition. They obtain some persistence related results based on Morse theory. It could be that there is a different way to perform the pairing in Section~\ref{chapter-pairing} using discrete Morse theory instead.

There are many potential applications to data analysis for this algorithm. In \cite{belchi_lung_2018}, persistent homology with a height filtration was utilized to detect disease in lungs. The full $3$D $\XPHT$ would provide a more complete picture of the problem. The $\XPHT_0$ part of this algorithm has already been used in \cite{HaruhisaOda} where it was benchmarked against other methods, and used to classify left and right lungs.

The question of applying the methods in here to $4$D manifolds naturally arises. There are a few obstructions to this. In the case of $4$D data, it usually stems from $3$D spatial data collected over time. This does not suggest it would have the structure of a compact manifold. Another obstruction arises from our methods that were $3$D specific in computing $H_1$ and the pairings. 

There are several steps where we chose our methods for the sake of obtaining a more straightforward way to check for correctness. For example, in the pairing section \ref{chapter-pairing}, we could have done the cohomology computation in a way similar to \cite{EdelsbrunnerOlsbock2020}. For future work, we could compare methodology to find an optimal algorithm for computing the $\XPHT$.

\bibliography{refs}

\end{document}